\documentclass[12pt]{article}
\usepackage{amsmath}
\usepackage{epsfig}
\usepackage{verbatim}
\usepackage{amssymb}
\usepackage{amsfonts}
\usepackage{amsbsy}
\usepackage{epsfig}

\newcommand{\nbyn}{{N\times N}}

\newcommand{\flip}{{{\|f\|}_{Lip}}}
\newcommand{\al}{{\alpha}}

\newcommand{\R}{{\mathbb R}}

\newcommand{\sigmaS}{{\sigma(S^{N^2-1})}}

\newtheorem{prop}{Proposition}[section]
\newtheorem{lem}[prop]{Lemma}

\newtheorem{remk}[prop]{Remark}
\newtheorem{coro}[prop]{Corollary}
\newtheorem{theo}[prop]{Theorem}

\begin{document}

\title {Concentration of the Spectral Measure for Large Random Matrices with Stable Entries}
\author{Christian Houdr\'e \thanks{Georgia Institute of Technology, School of Mathematics, Atlanta, Georgia, 30332-0160, houdre@math.gatech.edu} \and Hua Xu \thanks {Georgia Institute of
Technology, School of Mathematics, Atlanta, Georgia, 30332-0160,
xu@math.gatech.edu} }

\maketitle

\vspace{0.5cm}

\begin{abstract}
\noindent We derive concentration inequalities for functions of the
empirical measure of large random matrices with infinitely divisible
entries and, in particular, stable ones. We also give concentration
results for some other functionals of these random matrices, such as
the largest eigenvalue or the largest singular value.

\end{abstract}



\section{Introduction and Statements of Results:}

Large random matrices have recently attracted a lot of attention in
fields such as statistics, mathematical physics or combinatorics
(e.g., see Mehta [\ref{r28}], Bai and Silverstein [\ref{r38}],
Johnstone [\ref{r25}], Anderson, Guionnet and Zeitouni [\ref{r20}]).
For various classes of matrix ensembles, the asymptotic behavior of
the, properly centered and normalized, spectral measure or of the
largest eigenvalue is understood. Many of these results hold true
for matrices with independent entries satisfying some moment
conditions (Wigner [\ref{r19}], Tracy and Widom [\ref{r18}],
Soshnikov [\ref{r16}], Girko [\ref{r3}], Pastur [\ref{r13}], Bai
[\ref{r36}], G$\ddot{\text{o}}$tze and Tikhomirov [\ref{r39}]).

There is relatively little work outside the independent or finite
second moment assumptions. Let us mention Soshnikov [\ref{r15}] who,
using the method of determinants, studied the distribution of the
largest eigenvalue of Wigner matrices with entries having heavy
tails. (Recall that a real (or complex) Wigner matrix is a symmetric
(or Hermitian) matrix whose entries ${\bf M}_{i,i},{1\le i\le N}$,
and ${\bf M}_{i,j},{1\le i<j\le N}$, form two independent families
of iid (complex valued in the Hermitian case) random variables.) In
particular, (see [\ref{r15}]), for a properly normalized Wigner
matrix with entries belonging to the domain of attraction of an
$\al$-stable law, $\lim_{N\rightarrow
\infty}\mathbb{P}^N(\lambda_{max}\le x)=\exp{(-x^{-\al})}$ (here
$\lambda_{max}$ is the largest eigenvalue of such a normalized
matrix). Soshnikov and Fyodorov [\ref{r17}] further derived results
for the largest singular value of $K\times N$ rectangular random
matrices with independent Cauchy entries, showing that the largest
singular value of such a matrix is of order $K^2N^2$.

On another front, Guionnet and Zeitouni [\ref{r4}], gave
concentration results for functionals of the empirical spectral
measure for random matrices whose entries are independent and either
satisfy a Logarithmic Sobolev inequality or are compactly supported.
They obtained in that context, the subgaussian decay of the tails of
the empirical spectral measure when deviating from its mean (see
also Ledoux [\ref{r9}]). Our purpose in the present work is to deal
with matrices whose entries form a general infinitely divisible
vector, and in particular a stable one. We obtain concentration
results for functionals of the corresponding empirical spectral
measure, allowing for any type of light or heavy tails. The
methodologies developed here apply as well for the largest
eigenvalue or for the spectral radius of such random matrices.

Following the lead of Guionnet and Zeitouni [\ref{r4}], let us start
by setting our notation and framework.

 Let ${\mathcal M}_{\nbyn}(\mathbb C)$ be the set of $N\times N$ Hermitian matrices with complex entries,
 throughout, equipped with the Hilbert-Schmidt norm
 $$\|{\bf M}\|_{HS}=\sqrt{tr({\bf M}^*{\bf M})}=\sqrt{\sum_{i,j=1}^N|{\bf M}_{i,j}|^2}.$$
  Let $f$ be a real valued function on $\R$. The function $f$ can be
 viewed as mapping ${\mathcal M}_{\nbyn}(\mathbb{C})$ to ${\mathcal M}_{\nbyn}
 (\mathbb C)$. Indeed, for ${\bf M}=({\bf M}_{i,j})_{1\le i,j\le N}\in {\mathcal M}_{\nbyn}(\mathbb{C})$, so
 that ${\bf M}\!=\!{\bf U}{\bf D}{\bf U}^*$, where $\bf D$ is a diagonal matrix, with real entries $\lambda_1,...,\lambda_N$,
and $\bf U$ is a unitary
 matrix, set
 $$f({\bf M})={\bf U}f({\bf D}){\bf U}^*,\ f({\bf D})=\left(\begin{array}{clrr}     f(\lambda_1) & 0 & \cdots & 0 \\
     0 & f(\lambda_2) & \cdots & 0 \\       \vdots  & \vdots & \ddots  & \vdots \\ 0  & 0 & \cdots & f(\lambda_N) \\   \end{array}\right).$$
 Let $tr({\bf M})=\sum_{i=1}^N{\bf M}_{i,i}$ be the trace operator on ${\mathcal M}_{\nbyn}(\mathbb{C})$ and
 set also $$tr_N({\bf M})=\frac{1}{N}\sum_{i=1}^N{\bf M}_{i,i}.$$
 For a $N\times N$ random Hermitian
matrix with eigenvalues $\lambda_1,\lambda_2,...,\lambda_N$, let
$F_N(x)=\frac{1}{N}\sum_{i=1}^N{\mathbf 1}_{\{\lambda_i\le
 x\}}$ be the corresponding empirical spectral distribution function.
 As well known, if ${\bf M}$ is a $\nbyn$ Hermitian Wigner matrix with $\mathbb E[{\bf M}_{1,1}]=\mathbb E[{\bf M}_{1,2}]=0$,
 $\mathbb E[|{\bf M}_{1,2}|^2]=1$,
 and $\mathbb E[{\bf M}_{1,1}^2]<\infty$, the spectral
 measure of ${\bf M}/\sqrt N$ converges to the semicircle law: $\sigma(dx)=\sqrt{4-x^2}\mathbf{1}_{\{|x|\le
 2\}}dx/2\pi$ ([\ref{r20}]).

We study below the tail behavior of either the spectral measure or
the linear statistic of $f({\bf M})$ for classes of matrices ${\bf
M}$. Still following Guionnet and Zeitouni, we focus on a general
random matrix ${\bf {X_A}}$ given as follows:
 $${\bf {X_A}}=(({\bf {X_A}})_{i,j})_{1\le i,j\le N},\ {\bf {X_A}}={\bf {X_A^*}},\ ({\bf {X_A}})_{i,j}=\frac{1}{\sqrt N}A_{i,j}\omega_{i,j}, $$
 with $(\omega_{i,j})_{1\le i,j\le N}=(\omega_{i,j}^R+\sqrt{-1}\omega_{i,j}^I)_{1\le i,j\le
 N}$, $\omega_{i,j}=\overline{\omega_{j,i}}$,
 and where $\omega_{i,j}$, $1\le i\le j\le N$ is a complex valued random variable with law
$P_{i,j}=P_{i,j}^R+\sqrt {-1}P_{i,j}^I$, ${1\le i\le j\le N}$, with
$P_{i,i}^I=\delta_0$ (by the Hermite property). Moreover, the matrix
${\bf A}=({ A}_{i,j})_{1\le i,j\le N}$ is Hermitian with, in most
cases, non-random complex
 valued entries uniformly bounded, say, by $a$.

 Different choices for the
 entries of $\bf A$ allow to cover various
 types of ensembles. For instance, if $\omega_{i,j},{1\le i<j\le
 N}$, and $\omega_{i,i},{1\le i\le N}$, are iid $N(0,1)$ random variables,
taking $A_{i,i}=\sqrt 2$ and $A_{i,j}=1$, for $1\le i<j\le N$ gives
the GOE (Gaussian Orthogonal Ensemble). If
$\omega_{i,j}^R,\omega_{i,j}^I,{1\le i<j\le
 N}$, and $\omega_{i,i}^R,{1\le i\le N}$, are iid $N(0,1)$ random variables,
 taking $A_{i,i}=1$ and $A_{i,j}=1/\sqrt 2$, for $1\le i<j\le N$ gives the GUE (Gaussian Unitary
 Ensemble) (see [\ref{r28}]). Moreover, if $\omega_{i,j}^R,\omega_{i,j}^I,{1\le i<j\le
 N}$, and $\omega_{i,i}^R,{1\le i\le N}$, are two independent
 families of real valued random variables, taking $A_{i,j}=0$ for $|i-j|$ large and
 $A_{i,j}=1$ otherwise, gives band matrices. Proper choices of non-random $A_{i,j}$ also
make it possible to cover Wishart matrices, as seen in the later
part of this section. In certain instances, $A$ can also be chosen
to be random, like in the case of diluted matrices, in which case
$A_{i,j},{1\le i\le j\le N}$, are iid Bernoulli random variables
(see [\ref{r4}]).

On $\R^{N^2}$, let $\mathbb P^N$ be the joint law of the random
vector $X=({\omega_{i,i}^R},\omega_{i,j}^R,\omega_{i,j}^I)$, ${1\le
i<j\le N}$, where it is understood that the indices for
$\omega_{i,i}^R$ are $1\le i\le N$. Let $\mathbb E^N$ be the
corresponding expectation. Denote by $\hat{\mu}_A^N$ the empirical
spectral measure of the eigenvalues of ${\bf {X_A}}$, and further
note that
$$tr_Nf({\bf {X_A}})=\frac{1}{N}tr(f({\bf {X_A}}))=\int_ \mathbb R f(x)\hat{\mu}_{\bf A}^N(dx),$$ for any
bounded Borel function $f$. For a Lipschitz function
 $f: \mathbb{R}^d\rightarrow\mathbb{R}$, set
 $$\flip =\underset{x\ne y}{\sup}\frac{|f(x)-f(y)|}{\|x-y\|},$$
where throughout $\|\cdot\|$ is the Euclidean norm, and where we
write $f\in Lip(c)$ whenever $\flip\le c$.

Each element ${\bf M}$ of $\mathcal M_{N\times N}(\mathbb C)$ has a
unique collection of eigenvalues $\lambda=\lambda({\bf
M})=(\lambda_1,\cdots,\lambda_N)$ listed in non increasing order
according to multiplicity in the simplex
$$\mathcal S^N=\{\lambda_1\ge\cdots\ge \lambda_N: \lambda_i\in \R, 1\le i\le N\},$$
where throughout $\mathcal S^N$ is equipped with the Euclidian norm
$\|\lambda\|=\sqrt{\sum_{i=1}^N\lambda_i^2}.$ It is a classical
result sometimes called Lidskii's theorem ([\ref{r23}]), that the
map $\mathcal M_{N\times N}(\mathbb C)\rightarrow \mathcal S^N$
which associates to each Hermitian matrix its ordered list of real
eigenvalues is 1-Lipschitz ([\ref{r21}], [\ref{r10}]). For a matrix
${\bf {X_A}}$ under consideration with eigenvalues $\lambda({\bf
{X_A}})$, it is then clear that the map
$\varphi:({\omega_{i,i}^R},\omega_{i,j}^R,\omega_{i,j}^I)_{1\le
i<j\le N}\mapsto \lambda({\bf {X_A}})$ is Lipschitz, from $(\mathbb
R^{N^2},\|\cdot\|)$ to $(\mathcal S^N,\|\cdot\|)$, with Lipschitz
constant bounded by $a\sqrt{2/N}$. Moreover, for any real valued
Lipschitz function $F$ on $\mathcal S^N$ with Lipschitz constant
$\|F\|_{Lip}$, the map $F\circ\varphi$ is Lipschitz, from
$(\R^{N^2},\|\cdot\|)$ to $\R$, with Lipschitz constant at most
$a\|F\|_{Lip}\sqrt{2/N}$. Appropriate choices of $F$ ([\ref{r10}],
[\ref{r20}]) ensure that the maximal eigenvalue $\lambda_{max}({\bf
{X_A}})=\lambda_1({\bf {X_A}})$, the spectral radius $\rho({\bf
{X_A}})=\underset{1\le i\le N}{\max}|\lambda_i|$ and $tr_N(f({\bf
{X_A}}))$, where $f:\R\rightarrow\R$ is a Lipschitz function, are
themselves Lipschitz functions with Lipschitz constants at most
$a\sqrt {2/N}$, $a\sqrt {2/N}$ and $\sqrt 2a\flip/N$, respectively.
These observations (and our results) are also valid for the real
symmetric matrices, with proper modification of the Lipschitz
constants.

Next, Recall that $X$ is a $d$-dimensional infinitely divisible
random vector without Gaussian component, $X\sim ID(\beta,0,\nu)$,
if its characteristic function is given by,
\begin{align}\label{f1.1}
\varphi_X(t)&=\mathbb E e^{i\langle t,X\rangle} \nonumber\\
             &=\exp\left\{i\langle t,\beta\rangle+\int_{\mathbb{R}^d}\left(e^{i\langle t,u\rangle}-1-i\langle t,u\rangle\mathbf 1_{\|u\|\le1}\right)\nu(du)\right\},
\end{align}
where $t,\beta\in \mathbb{R}^d$ and $\nu\not \equiv 0$ (the
$L\acute{e}vy\ measure$) is a positive measure on
$\mathcal{B}(\mathbb{R}^d)$, the Borel $\sigma$-field of $\R^d$,
without atom at the origin, and such that
$\int_{\mathbb{R}^d}(1\wedge\|u\|^2)\nu(du)<+\infty$. The vector $X$
has independent components if and only if its L\'evy measure $\nu$
is supported on the axes of $\R^d$ and is thus of the form:
\begin{equation}
\nu(dx_1,\dots,dx_d)\!=\!\!\sum_{k=1}^{d}\delta_0(dx_1)\dots
\delta_0(dx_{k-1})\tilde{\nu}_k(dx_k)\delta_0(dx_{k+1})
\dots\delta_0(dx_{d}),
\end{equation}
for some one-dimensional L\'evy measures $\tilde{\nu}_k$. Moreover,
the $\tilde{\nu}_k$ are the same for all $k=1,\dots,d$, if and only
if $X$ has identically distributed components.

The following proposition gives an estimate on any median (or the
mean, if it exists) of a Lipschitz function of an infinitely
divisible vector $X$. It is used in most of the results presented in
this paper. The first part is a consequence of Theorem 1 in
[\ref{r6}], while the proof of the second part can be obtained as in
[\ref{r6}].

\begin{prop}\label{prop1.1}
Let $X=(\omega_{i,i}^R,\omega_{i,j}^R,\omega_{i,j}^I)_{1\le i<j\le
N}\sim ID(\beta,0,\nu)$ in $\R^{N^2}$. Let $V^2(x)=\int_{\|u\|\le
x}{{\|u\|}^2\nu(du)}$, $\bar{\nu}(x)=\int_{\|u\|>x}\nu(du)$, and for
any $\gamma>0$, let $p_\gamma=\inf\big\{x>0:0<V^2(x)/x^2\le
\gamma\big\}$. Let $f\in Lip(1)$, then for any $\gamma$ such that
$\bar{\nu}(p_\gamma)\le 1/4$,
\begin{itemize}
\item[(i)] any median $m(f(X))$ of $f(X)$ satisfies
$$
|m(f(X))-f(0)|\le
G_1(\gamma):=p_\gamma\Big(\sqrt{\gamma}+3k_\gamma(1/4)\Big)+E_\gamma,
$$

\item[(ii)] the mean $\mathbb E^N[f(X)]$ of $f(X)$, if it exists, satisfies
$$
|\mathbb E^N[f(X)]-f(0)|\le
G_2(\gamma):=p_\gamma\Big(\sqrt{\gamma}+k_\gamma(1/4)\Big)+E_\gamma,
$$

\end{itemize}
\noindent where $k_\gamma(x)$, $x>0$, is the solution, in y, of the
equation
$$y-\left(y+\gamma\right)\ln\left(1+\frac{y}{\gamma}\right)=\ln x,$$
\noindent and where
\begin{align}
E_\gamma=\Bigg(\sum_{k=1}^{N^2}\Big(\langle e_k,b\rangle
-\int_{p_\gamma<\|y\|\le1}\langle e_k,y\rangle
\nu(dy)+\int_{1<\|y\|\le p_\gamma}\langle e_k,y\rangle
\nu(dy)\Big)^2\Bigg)^{1/2},
\end{align}

\noindent with $e_1,e_2,\dots,e_{N^2}$ being the canonical basis of
$\R^{N^2}$.

\end{prop}

Our first result deals with the spectral measure of a Hermitian
matrix whose entries on and above the diagonal form an infinitely
divisible random vector with finite exponential moments. Below, for
any $b>0$, $c>0$, let
$${Lip}_ b(c)=\Big\{f:\R\rightarrow
\R\ :\flip\le c,\ \|f\|_\infty\le b\Big\},$$

\noindent while for a fixed compact set $\mathcal K\subset \R$, with
diameter $|\mathcal K|=\underset{x,y\in\mathcal K}{\sup}|x-y|$, let
$${Lip_{\mathcal K}(c)}:=\{f:\R\rightarrow
\R\ :\flip\le c,\ supp(f)\subset\mathcal K\},$$ where $supp(f)$ is
the support of $f$.

\begin{theo}\label{theo1.1}
Let $X=(\omega_{i,i}^R,\omega_{i,j}^R,\omega_{i,j}^I)_{1\le i<j\le
N}$ be a random vector with joint law $\mathbb{P}^N\sim
ID(\beta,0,\nu)$ such that $\mathbb{E}^N [e^{t\|X\|}]<+\infty$, for
some $t>0$. Let $T=\sup\{t\ge0:\mathbb{E}^N
\big[e^{t\|X\|}\big]<+\infty\}$ and let $h^{-1}$
 be the inverse of $$h(s)=\int_{\R^{N^2}}\|u\| \big(e^{s\|u\|}-1\big)\nu(du),\ \ \ \ 0<s<T.$$
\begin{itemize}
\item[(i)] For any compact set $\mathcal K\subset \mathbb R$,
\begin{align}\label{f1.3}
\mathbb{P}^N \Bigl(\underset{f\in {Lip_{\mathcal K}(1)}}{\sup}|tr_N(f({\bf {X_A}}&))-\mathbb E^N\left[tr_N(f({\bf {X_A}}))\right]|\ge\delta\Bigr) \nonumber \\
&\le \frac{8|\mathcal K|}{\delta}\exp
\Bigg\{-\int_0^{\frac{N\delta^2}{8\sqrt2a|\mathcal
K|}}h^{-1}(s)ds\Bigg\},
\end{align}
for all $\delta>0$ such that $\delta^2< 8\sqrt2a|\mathcal
K|h\left(T^{-}\right)/N$.

\item[(ii)]
\begin{align}\label{f1.4}
\mathbb{P}^N \Bigl(\underset{f\in Lip_b(1)}{\sup}|tr_N(f({\bf
{X_A}}&))-\mathbb E^N\left[tr_N(f({\bf
{X_A}}))\right]|\ge\delta\Bigr)\nonumber\\
&\le\frac{C(\delta,b)}{\delta}\exp\!\Bigg\{-\int_0^{\frac{N\delta^2}{\sqrt2aC(\delta,b)}}h^{-1}(s)ds\Bigg\},
\end{align}
\noindent for all $\delta>0$ such that $\delta^2\le
\sqrt2aC(\delta,b)h(T^-)/N$, where
$$C(\delta,b)=C\bigg(\frac{\sqrt 2a}{\sqrt N}\Big(G_2(\gamma)+h(t_0)\Big)+b\bigg),$$
\noindent with $G_2(\gamma)$ as in Proposition \ref{prop1.1}, $C$ a
universal constant, and with $t_0$ the solution, in $t$, of
$th(t)-\int_0^th(s)ds-\ln(12b/\delta)=0$.

\end{itemize}

\end{theo}

\begin{remk}
\begin{itemize}
\item[(i)] The order of $C(\delta,b)$ in part (ii) can be made more specific. Indeed, it will be clear from the proof of the theorem (see
(\ref{f2.26})), that for any $0<t^*\le T$ fixed,
$$C(\delta,b)\le C\bigg(\frac{\sqrt 2a}{\sqrt N}\Big(\frac{\ln\frac{12b}{\delta}}{t^*}+\frac{\int_0^{t^*}h(s)ds}{t^*}+G_2(\gamma)\Big)\bigg).$$

\item[(ii)] As seen from the proof (see (\ref{f2.12})), in the statement of the above theorem, $G_2(\gamma)$ can be replaced by
$\mathbb E^N\big[\|X\|\big]$. Now $\mathbb E^N\big[\|X\|\big]$ is of
order $N$ since,
\begin{equation}\label{f1.12}
N \underset{j=1,2,\dots,N^2}{\min}\mathbb
E^N\big[|X_j|\big]\le\mathbb E^N\big[\|X\|\big]\le N
\underset{j=1,2,\dots,N^2}{\max}\mathbb E^N\big[X_j^2\big],
\end{equation}
\noindent where the ${X_j}$, $j=1,2,\dots,N^2$ are the components of
$X$. Actually, an estimate more precise than (\ref{f1.12}) is given
by a result of Marcus and Rosi\'nski [\ref{r37}] which asserts that
if $\mathbb E[X]=0$, then
$$\frac{1}{4}x_0\le \mathbb E\big[\|X\|\big]\le \frac{17}{8}x_0,$$
\noindent where $x_0$ is the solution of the equation:
\begin{equation}\label{e2.1}
\frac{V^2(x)}{x^2}+\frac{M(x)}{x}=1,
\end{equation}
\noindent where $V^2(x)$ is as before, while $ M(x)=\int_{\|u\|\ge
x}{{\|u\|}\nu(du)}$, $x>0$.

\item[(iii)] As usual, one can easily pass from the mean $\mathbb E^N[tr_N(f)]$
to any median $m(tr_N(f))$ in either (\ref{f1.3}) or (\ref{f1.4}). Indeed, for any
$0\le\delta\le 2b$, if $$\underset{f\in
Lip_b(1)}{\sup}|tr_N(f)-m(tr_N(f))|\ge\delta,$$ there exist a
function $f\in Lip_b(1)$ and a median $m(tr_N(f))$ of $tr_N(f)$,
such that either $tr_N(f)-m(tr_N(f))\ge\delta$ or
$tr_N(f)-m(tr_N(f))\le-\delta$. Without loss of generality assuming
the former, otherwise dealing with the latter with $-f$, consider
the function $g(y)=\min\left(d(y, A),\delta\right)/2$,
$y\in\R^{N^2}$, where $A=\left\{tr_N(f)\le m(tr_N(f)\right\}$.
Clearly $g\in Lip_b(1)$, $\mathbb E^N[tr_N(g)]\le \delta/4$, and
therefore $tr_N(g)-\mathbb E^N\left[tr_N(g)\right]\ge\delta/4$,
which indicates that
$$\underset{g\in Lip_b(1)}{\sup}\!\big|tr_N(g)-\mathbb
E^N\left[tr_N(g)\right]\big|\ge\frac{\delta}{4}.$$ \noindent Hence,
\begin{align}
\mathbb P^N\bigg(\underset{f\in
Lip_b(1)}{\sup}\big|tr_N&(f)-m(tr_N(f))\big|\ge\delta\bigg)\nonumber\\
&\le\mathbb P^N\bigg(\underset{g\in
Lip_b(1)}{\sup}\big|tr_N(g)-\mathbb
E^N\left[tr_N(g)\right]\big|\ge\frac{\delta}{4}\bigg).
\end{align}
\end{itemize}
\end{remk}

Next, recall (see [\ref{r26}], [\ref{r10}]) that the Wasserstein
distance between any two probability measures $\mu_1$ and $\mu_2$ on
$\R$ is defined by
\begin{equation}
\begin{split}
d_W(\mu_1,\mu_2)&=\underset{f\in Lip_b(1)}{\sup}\Big|\int_\R
fd\mu_1-\int_\R fd\mu_2\Big|.
\end{split}
\end{equation}
\noindent \noindent Hence, Theorem \ref{theo1.1} actually gives a
concentration result, with respect to the Wasserstein distance, for
the empirical spectral measure $\hat{\mu}_{\bf A}^N$, when it
deviates from its mean $\mathbb E^N[\hat{\mu}_{\bf A}^N]$.

As in [\ref{r4}], we can also obtain a concentration result for the
distance between any particular probability measure and the
empirical spectral measure.

\begin{prop}\label{theo1.3}
Let $X=(\omega_{i,i}^R,\omega_{i,j}^R,\omega_{i,j}^I)_{1\le i<j\le
N}$ be a random vector with joint law $\mathbb{P}^N\sim
ID(\beta,0,\nu)$ such that $\mathbb{E}^N
\big[e^{t\|X\|}\big]<+\infty$, for some $t>0$. Let
$T=\sup\{t>0:\mathbb{E}^N \big[e^{t\|X\|}\big]<+\infty\}$ and let
$h^{-1}$ be the inverse of $h(s)=\int_{\R^{N^2}}{\|u\|} (e^{s
\|u\|}-1)\nu(du),\ 0<s<T$. Then, for any probability measure $\mu$,
\begin{equation}
\mathbb{P}^N\left(d_W(\hat{\mu}_A^N, \mu)-\mathbb
E^N[d_W(\hat{\mu}_A^N, \mu)]\ge \delta\right)\le
\exp\bigg\{-\int_0^{\frac{N\delta}{\sqrt2a}}h^{-1}(s)ds\bigg\},
\end{equation}

\noindent for all $0<\delta< \sqrt2a h\left(T^{-}\right)/N$.
\end{prop}

Of particular importance is the case of an infinitely divisible
vector having boundedly supported L\'evy measure. We then have:

\begin{coro}\label{coro1.4}
Let $X=(\omega_{i,i}^R,\omega_{i,j}^R,\omega_{i,j}^I)_{1\le i<j\le
N}$ be a random vector with joint law $\mathbb{P}^N\sim
ID(\beta,0,\nu)$ such that $\nu$ has bounded support. Let
$R=\inf\{r>0:\nu(x:\|x\|>r)=0\}$, let
$V^2\big(=V^2(R)\big)=\int_{\R^{N^2}}\|u\|^2\nu(du)$, and let
$$\ell(x)=(1+x)\ln(1+x)-x,$$
\noindent $x>0$.
\begin{itemize}
\item[(i)] For any $\delta>0$,
\begin{align}\label{f1.13}
\mathbb{P}^N \Bigl(\underset{f\in Lip_b(1)}{\sup}|tr_N&(f({\bf {X_A}}))-\mathbb E^N[tr_N(f({\bf {X_A}}))]|\ge\delta\Bigr) \nonumber \\
\noindent &\le \frac{C(\delta,b)}{\delta}\exp
\Bigg\{-\frac{V^2}{R^2}\ell
\bigg(\frac{NR\delta^2}{\sqrt2aC(\delta,b)V^2}\biggr)\Bigg\},
\end{align}
\noindent where
$$C(\delta,b)=C\bigg(\frac{\sqrt 2a}{\sqrt N}\Big(G_2(\gamma)+\frac{V^2}{R}\big(e^{t_0R}-1\big)\Big)+b\bigg),
$$
\noindent with $G_2(\gamma)$ as in Proposition \ref{prop1.1}, $C$ a
universal constant, and $t_0$ the solution, in $t$, of
$$\frac{V^2}{R^2}\Big(tRe^{tR}-e^{tR}+1\Big)=\ln\frac{12b}{\delta}.$$

\item[(ii)] For any probability measure $\mu$ on $\R$, and any $\delta>0$,
\begin{align}
\mathbb{P}^N&\bigl(d_W(\hat{\mu}_A^N, \mu)-\mathbb
E^N[d_W(\hat{\mu}_A^N, \mu)]\ge \delta\bigr)\nonumber \\
&\le \exp \Bigg\{\frac{N\delta}{\sqrt
2aR}-\Biggl(\frac{N\delta}{\sqrt
2aR}+\frac{V^2}{R^2}\Biggr)\ln\Biggl(1+\frac{NR\delta^2}{\sqrt
2aV^2}\Biggr)\Bigg\}.
\end{align}

\end{itemize}

\end{coro}

\begin{remk}
\begin{itemize}
\item[(i)] As in Theorem \ref{theo1.1}, the dependency of $C(\delta,b)$ in
$\delta$ and $b$ can be made more precise. A key step in the proof
of (\ref{f1.13}) is to choose $\tau$ such that
$$\mathbb E^N[tr_N(\mathbf 1_{\{|{\bf{X_A}}|\ge\tau\}})]\le
\delta/12b,$$ \noindent and then $C(\delta,b)$ is determined by
$\tau$. Minimizing, in $t$, the right hand side of (\ref{f2.12}),
leads to the following estimate
$$\mathbb E^N[tr_N(\mathbf 1_{\{|{\bf
{X_A}}|\ge\tau\}})]\le\exp\Bigg\{-\frac{V^2}{R^2}\ell\Bigg(\frac{R\Big(\frac{\sqrt
N}{\sqrt 2a}\tau-G_2(\gamma)\Big)}{V^2}\Bigg)\Bigg\},$$ \noindent
where $\ell(x)=(1+x)\ln(1+x)-x$. For $x\ge1$, $2\ell(x)\ge x\ln x$.
Hence one can choose $\tau$ to be the solution, in $x$, of the
equation
$$\frac{x}{R}\ln\frac{xR}{V^2}=2\ln\frac{12b}{\delta}.$$
\noindent It then follows that $C(\delta,b)$ can be taken to be
$$C\bigg(\frac{\sqrt 2a}{\sqrt N}\Big(G_2(\gamma)+\tau\Big)+b\bigg).$$
\end{itemize}
\end{remk}

Outside of the finite exponential moment assumption, an interesting
class of random matrices with infinitely divisible entries are the
ones with stable entries, which we now analyze.

Recall that $X$ in $\R^d$ is $\alpha$-stable, $(0<\alpha<2)$, if its L\'evy
measure $\nu$ is given, for any Borel set $B\in \mathcal B(\mathbb
R^d)$, by
\begin{equation}
\label{f1.6}
\nu(B)=\int_{S^{d-1}}\sigma(d\xi)\int_0^{+\infty}\mathbf
1_B(r\xi)\frac{dr}{r^{1+\alpha}},
\end{equation}
where $\sigma$, the spherical component of the L\'evy measure, is a
finite positive measure on $S^{d-1}$, the unit sphere of $\mathbb
R^d$. Since the expected value of the spectral measure of a matrix
with $\alpha$-stable entries might not exist, we look at the
deviation from a median. Here is a sample result.

\begin{theo}\label{theo1.5}
Let $X=(\omega_{i,i}^R,\omega_{i,j}^R,\omega_{i,j}^I)_{1\le i<j\le
N}$ be an $\alpha$-stable, $0<\al<2$, random vector in
$\R^{N^2}\!\!$ with L\'evy measure $\nu$ given by (\ref{f1.6}).

\begin{itemize}
\item[(i)]
Let $f\in Lip(1)$, and let $m(tr_N(f({\bf {X_A}})))$ be any median
of $tr_N(f({\bf {X_A}}))$. Then,
\begin{equation}\label{f1.9}
\mathbb{P}^N\!\bigl(tr_N(f({\bf {X_A}}\!))-m(tr_N(f({\bf
{X_A}}\!)))\!\ge\!\delta\bigr)\!\le\! C(\al)(\sqrt
2a)^\al\frac{\sigma(S^{N^2-1})}{N^\al\delta^\al},
\end{equation}

\noindent whenever $\delta N>\sqrt
2a\left[2\sigma(S^{N^2-1})C(\al)\right]^{1/\al}$, and where
$C(\al)=4^\al(2-\al+e\al)/\al(2-\al)$.

\item[(ii)]
Let $\lambda_{max}({\bf {X_A}})$ be the largest eigenvalue of ${\bf
{X_A}}$, and let $m(\lambda_{max}({\bf {X_A}}))$ be any median of
$\lambda_{max}({\bf {X_A}})$, then
\begin{equation}\label{f1.10}
\mathbb{P}^N\bigl(\lambda_{max}({\bf {X_A}})-m(\lambda_{max}({\bf
{X_A}}))\ge\delta\bigr)\le C(\al)(\sqrt
2a)^\al\frac{\sigma(S^{N^2-1})}{N^{\al/2}\delta^\al},
\end{equation}

\noindent whenever $\delta \sqrt N>\sqrt
2a\left[2\sigma(S^{N^2-1})C(\al)\right]^{1/\al}$, and where
$C(\al)=4^\al(2-\al+e\al)/\al(2-\al)$.

\end{itemize}
\end{theo}

\begin{remk}

Let ${\bf M}$ be a Wigner matrix whose entries ${\bf M}_{i,i}$,
${\!1\!\!\le\!\! i\!\!\le\! \! N}$, $\!{\bf M}_{i,j}^R$,
${\!1\!\!\le\!\! i\!\!<\!\!j\!\!\le\!\! N}$, and ${\bf M}_{i,j}^I$,
${1\!\le\! i\!<\!j\!\le\! N}$, are iid random variables, such that
the distribution of $|{\bf M}_{1,1}|$ belongs to the domain of
attraction of an $\al$-stable distribution, i.e., for any
$\delta>0$,
$$\mathbb P(|{\bf M}_{1,1}|>\delta)=\frac{L(\delta)}{\delta^\al},$$
\noindent for some slowly varying positive function $L$ such that
$\underset {\delta\rightarrow\infty}{\lim}L(t\delta)/L(\delta)=1$,
for all $t>0$. Soshnikov [\ref{r15}] showed that, for any
$\delta>0$,
$$\lim_{N\rightarrow \infty}\mathbb P^N(\lambda_{max}(b_N^{-1}{\bf M})\ge \delta)=1-\exp(-\delta^{-\al}),$$
\noindent where $b_N$ is a normalizing factor such that $\underset
{N\rightarrow\infty}{\lim}N^2L(b_N)/b_N^\al=2$ and where
$\lambda_{max}(b_N^{-1}{\bf M})$ is the largest eigenvalue of
$b_N^{-1}{\bf M}$. In fact $\underset
{N\rightarrow\infty}{\lim}N^{\frac{2}{\al}-\epsilon}/b_N=0$ and
$\underset {N\rightarrow\infty}{\lim}
b_N/N^{\frac{2}{\al}+\epsilon}=0$, for any $\epsilon>0$. As stated
in [\ref{r7}], when the random vector $X$ is in the domain of
attraction of an $\al$-stable distribution, concentration
inequalities similar to (\ref{f1.9}) or (\ref{f1.10}) can be
obtained for general Lipschitz function. In particular, if the
L\'evy measure of $X$ is given by
 \begin{equation}
 \label{f1.7}
 \nu(B)=\int_{S^{N^2-1}}\sigma(d\xi)\int_0^{+\infty}\mathbf
 1_B(r\xi)\frac{L(r)dr}{r^{1+\alpha}},
 \end{equation}
for some slowly varying function $L$ on $[0,+\infty)$, and if we
still choose the normalizing factor $b_N$ such that
$\lim_{N\rightarrow\infty}\sigmaS L(b_N)/b_N^\al$ is constant, then,

\begin{align}
\mathbb P^N\bigl(\lambda_{max}(b_N^{-1}{\bf
M})-m(\lambda_{max}&(b_N^{-1}{\bf M}))\ge \delta\bigr) \nonumber \\
&\le \frac{C(\al)\sigmaS{2}^{\al/2}
}{b_N^\al}\frac{L\Big(b_N\frac{\delta}{\sqrt2}\Big)}{\delta^\al},
\end{align}
\noindent whenever
$${(\delta b_N)}^{\al}\ge 2^{1+\al/2}C(\al)\sigmaS L\big(b_N\delta/\sqrt2\big).$$
\noindent Now, recall that for an $N^2$ dimensional vector with iid
entries, $\sigma(S^{N^2-1})=N^2({\hat\sigma}(1)+{\hat\sigma}(-1))$,
where $\hat\sigma(1)$ is short for $\sigma(1,0,\dots,0)$ and
similarly for $\hat\sigma(-1)$. Thus, for fixed $N$, our result
gives the correct order of the upper bound for large values of
$\delta$, since for $\delta>1$,
$$\frac{e-1}{e\delta^\al}\le 1-e^{-\delta^{-\al}}\le \frac{1}{\delta^\al}.$$
\noindent Moreover, in the stable case, $L(\delta)$ becomes
constant, and $b_N=N^{2/\al}$. Since $\lambda_{max}(N^{-2/\al}{\bf
M})$ is a Lipschitz function of the entries of the matrix ${\bf M}$
with Lipschitz constant at most $\sqrt 2N^{-2/\al}$, for any median
$m(\lambda_{max}(N^{-2/\al}{\bf M}))$ of
$\lambda_{max}(N^{-2/\al}{\bf M})$, we have,
\begin{equation}\label{f1.17}
\mathbb P^N\!\bigl(\lambda_{max}(N^{-\frac{2}{\al}}{\bf
M})-m(\lambda_{max}(N^{-\frac{2}{\al}}{\bf M}))\!\ge\!
\delta\bigr)\!\le\!
C(\al)\frac{\big({\hat\sigma}(1)\!+\!{\hat\sigma}(-1)\big)}{2^{\al/2}}\frac{1}{\delta
^\al},
\end{equation}
 \noindent whenever $\delta\ge
\left[2C(\al)\big({\hat\sigma}(1)+{\hat\sigma}(-1)\big)\right]^{1/\al}$.
Furthermore, using Theorem 1 in [\ref{r6}], it is not difficult to
see that $m(\lambda_{max}(N^{-2/\al}{\bf M}))$ can be upper and
lower bounded independently of $N$. Finally, an argument as in
Remark \ref{remk1.12} below will give a lower bound on
$\lambda_{max}(N^{-2/\al}{\bf M})$ of the same order as
(\ref{f1.17}).
\end{remk}

\vspace{.15in} The following proposition will give an estimate on
any median of a Lipschitz function of $X$, where $X$ is a stable
vector. It is the version of Proposition \ref{prop1.1} for
$\al$-stable vectors.

\begin{prop}\label{prop2.5}
Let $X=(\omega_{i,i}^R,\omega_{i,j}^R,\omega_{i,j}^I)_{1\le i<j\le
N}$ be an $\alpha$-stable, $0<\al<2$, random vector in $\R^{N^2}$
with L\'evy measure $\nu$ given by (\ref{f1.6}). Let $f\in Lip(1)$,
then
\begin{itemize}
\item[(i)] any median $m(f(X))$ of $f(X)$ satisfies
\begin{align}
&|m(f(X))-f(0)| \nonumber\\
&\!\le\!
J_1(\al):=\!\left(\frac{\sigma(S^{N^2-1})}{4\alpha}\right)^{1/\alpha}\!\!\left(\sqrt{\frac{\al}{4(2-\al)}}+3k_{\frac{\al}{4(2-\al)}}(1/4)\right)\!+\!E,
\end{align}
\item[(ii)] the mean $\mathbb E^N[f(X)]$ of $f(X)$, if it exists,
satisfies
\begin{align}
&|\mathbb E^N[f(X)]-f(0)| \nonumber\\
&\!\le\!
J_2(\al):=\!\left(\frac{\sigma(S^{N^2-1})}{4\alpha}\right)^{1/\alpha}\!\!\left(\sqrt{\frac{\al}{4(2-\al)}}+k_{\frac{\al}{4(2-\al)}}(1/4)\right)\!+\!E,
\end{align}
\end{itemize}

\noindent where $k_{\al/4(2-\al)}(x)$, $x>0$, is the solution, in y,
of the equation
$$y-\left(y+\frac{\al}{4(2-\al)}\right)\ln\left(1+\frac{4(2-\al)y}{\al}\right)=\ln x,$$
\noindent and where
\begin{align}
E=\Bigg(\sum_{k=1}^{N^2}\Big(\langle e_k,b\rangle
-&\int_{\big(\frac{4\sigma(S^{N^2-1})}{\alpha}\big)^{1/\alpha}<\|y\|\le1}\langle
e_k,y\rangle \nu(dy)\nonumber\\
&+\int_{1<\|y\|\le\big(\frac{4\sigma(S^{N^2-1})}{\alpha}\big)^{1/\alpha}}\langle
e_k,y\rangle \nu(dy)\Big)^2\Bigg)^{1/2},
\end{align}

\noindent with $e_1,e_2,\dots,e_{N^2}$ being the canonical basis of
$\R^{N^2}$.

\end{prop}

\begin{remk}
\begin{itemize}

\item[(i)] When the components of $X$ are independent, a direct computation
shows that, up to a constant, $E$ in both $J_1(\al)$ and $J_2(\al)$
is dominated by
$\left(\frac{\sigma(S^{N^2-1})}{4\alpha}\right)^{1/\alpha}$, as
$N\rightarrow\infty$.
\end{itemize}
\end{remk}

In complete similarity to the finite exponential moments case, we
can obtain concentration results for the spectral measure of
matrices with $\al$-stable entries.

\begin{theo}\label{theo1.7}
Let $X=(\omega_{i,i}^R,\omega_{i,j}^R,\omega_{i,j}^I)_{1\le i<j\le
N}$ be an $\alpha$-stable, $0\!<\!\al\!<\!2$, random vector in
$\R^{N^2}\!\!$ with L\'evy measure $\nu$ given by (\ref{f1.6}).

\begin{itemize}
\item[(i)] Then,
\begin{align}
\mathbb{P}^N\biggl(\underset{f\in Lip_b(1)}{\sup}\big|tr_N(f({\bf
{X_A}}))-&\mathbb E^N[tr_N(f({\bf
{X_A}}))]\big|\ge\delta\biggr)\nonumber\\
&\le C(\delta,b,\al)\frac{a^{\al}\sigmaS}{N^\al\delta^{\al}}\wedge1,
\end{align}
\noindent where
$$C(\delta,b,\al)=\bigg(C_1(\al)\bigg(\frac{\sqrt2a}{\sqrt
N}\bigg)^{1+\al}\bigg(\frac{J_1(\al)\!+\!1}{\delta}+b\bigg)^{1+\al}+C_2(\al)\bigg),$$
\noindent with $C_1(\al)$ and $C_2(\al)$ constants depending only on
$\al$, and with $J_1(\al)$ as in Proposition \ref{prop2.5}.

\item[(ii)] For any probability measure
$\mu$,
\begin{equation}\label{f1.11}
\mathbb{P}^N\!\!\left(d_W(\hat{\mu}_A^N,
\mu)\!-\!m(d_W(\hat{\mu}_A^N, \mu))\!\ge\!\delta\right)\!\le\!
C(\al)(\sqrt 2a)^{\al}\frac{\sigma(S^{N^2-1})}{N^\al\delta^\al},
\end{equation}

\noindent whenever $\delta N\ge \sqrt
2a\left[2\sigma(S^{N^2-1})C(\al)\right]^{1/\al}$ and where
$C(\al)=4^\al(2-\al+e\al)/\al(2-\al)$.

\end{itemize}
\end{theo}

 \vspace{.14in}

It is also possible to obtain concentration results for smaller
values of $\delta$. The lower and intermediate range for the stable
deviation obtained in [\ref{r2}] provide the appropriate tools to
achieve the following result. We refer to [\ref{r2}] for complete
arguments, and only provide below a sample result.

\vspace{.15 in}
\begin{theo}\label{theo1.8}
Let $X=(\omega_{i,i}^R,\omega_{i,j}^R,\omega_{i,j}^I)_{1\le i<j\le
N}$ be an $\alpha$-stable, $1\!<\!\al\!<\!2$, random vector in
$\R^{N^2}\!\!$ with L\'evy measure $\nu$ given by (\ref{f1.6}). For
any $\epsilon>0$, there exists $\eta(\epsilon)$, and constants
$D_1=D_1\big(\al,a,N,\sigmaS\big)$ and
$D_2=D_2\big(\al,a,N,\sigmaS\big)$, such that for all
$0<\delta<\eta(\epsilon)$,
\begin{align}\label{f1.24}
\mathbb{P}^N\bigl(\underset{f\in Lip_b(1)}{\sup}\big|tr_N(f({\bf
{X_A}}))-\mathbb E^N[tr_N&(f({\bf
{X_A}}))]\big|\ge\delta\bigr)\nonumber \\
&\le(1+\epsilon)\frac{D_1}{\delta^{\frac{\al+1}{\al}}}\exp\Bigl(-D_2\delta^{\frac{2\al+1}{\al-1}}\Bigr).
\end{align}
\end{theo}
\vspace{.15 in}

\begin{remk}\label{remk1.13}
\begin{itemize}
\item[(i)] In (\ref{f1.9}), (\ref{f1.10}) or (\ref{f1.11}), the
constant $C(\al)$ is not of the right order as $\al\rightarrow 2$.
It is, however, a simple matter to adapt Theorem 2 of [\ref{r7}] to
obtain, at the price of worsening the range of validity of the
concentration inequalities, the right order in the constants as
$\al\rightarrow 2$.

\item[(ii)] Let us now provide some estimation of $D_1$ and $D_2$, which are
needed for comparison with the GUE results of [\ref{r4}] (see (iii)
below). Let $C(\al)=2^\al(e\al+2-\al)/(2(2-\al))$,
$K(\al)=\max\Big\{2^\al/(\al-1),C(\al)\Big\}$,
$L(\al)=\big((\al-1)/\al\big)^{\al/(\al-1)}(2-\al)/10$ and let
\begin{align}\label{f1.26}
D^*\!=\!2\bigg(\frac{\sqrt2a}{\sqrt
N}\bigg)^{\frac{2\al-1}{\al}}&\!\!\bigg(12\frac{C(\al)}{K(\al)}\bigg)^{\frac{1}{\al}}\!J_2(\al)b^{\frac{1}{\al}}+2\bigg(\frac{\sqrt2a}{\sqrt
N}\bigg)^{\frac{\al-1}{\al}}\!\!\bigg(12\frac{C(\al)}{K(\al)}\bigg)^{\frac{1}{\al}}\!b^{\frac{\al+1}{\al}}\nonumber\\
&+\frac{2\sqrt2a}{\sqrt
N}\Big(12C(\al)\sigmaS\Big)^{\frac{1}{\al}}b^{\frac{1}{\al}}.
\end{align}
\noindent As shown in the proof of the theorem, $ D_1=24D^*$, while
$$
D_2=\frac{L(\al)}{\Big(\sigmaS\Big)^{\frac{1}{\al-1}}}\bigg(\frac{N}{\sqrt2a}\bigg)^{\frac{\al}{\al-1}}
\frac{1}{\big(72D^*\big)^{\frac{\al}{\al-1}}}. $$ \noindent Thus, as
$N\rightarrow+\infty$, $D_1$ is of order
$N^{-1/2}\bigg(\sigmaS\bigg)^{1/\al}$, while $D_2$ is of order
$N^{3\al/(2\al-2)}\bigg(\sigmaS\bigg)^{2/(1-\al)}$.
\item[(iii)] Guionnet and Zeitouni [\ref{r4}], obtained concentration
results for the spectral measure of matrices with independent entries, which are either
compactly supported or satisfy a logarithmic Sobolev inequality. In
particular for the elements of the GUE, their upper bound of
concentration for the spectral measure is
\begin{equation}\label{f1.25}
\frac{C_1+b^{3/2}}{\delta^{3/2}}\exp\bigg\{-\frac{C_2}{8ca^2}N^2\frac{\delta^5}{(C_1+b^{3/2})^2}\bigg\},
\end{equation}
\noindent where $C_1$ and $C_2$ are universal constants. In Theorem
\ref{theo1.8}, the order, in $b$, of $D_1$ is at most
$b^{\al+1/\al}$, while that of $D_2$ is at least
$b^{-(\al+1)/(\al-1)}$. This order is thus consistent with the one
in (\ref{f1.25}), as $\al$ is close to 2. Taking into account part
(ii) above, the order of the constants in (\ref{f1.24}) are correct
when $\al\rightarrow 2$. Following [\ref{r2}] (see also Remark 4 in
[\ref{rM}]), we can recover a suboptimal Gaussian result by
considering a particular stable random vector $X^{(\al)}$ and
letting $\al\rightarrow 2$. Toward this end, let $X^{(\al)}$ be the
stable random vector whose L\'evy measure has for spherical
component $\sigma$, the uniform measure with total mass
$\sigmaS=N^2(2-\al)$. As $\al$ converges to 2, $X^{(\al)}$ converges
in distribution to a standard normal random vector. Also, as
$\al\rightarrow 2$, the range in $\delta$ in Theorem \ref{theo1.8}
becomes $(0,+\infty)$ while the constants in the concentration bound
do converge. Thus, the right hand side of (\ref{f1.24}) becomes
$$\frac{D_1}{\delta^{3/2}}\exp\bigg\{-D_2\delta^5\bigg\},$$
\noindent which is of the same order, in $\delta$, as (\ref{f1.25}).
However our order in $N$ is suboptimal.

\item[(iv)] In the proof of Theorem \ref{theo1.8}, the desired estimate in
(\ref{f2.54}) is achieved through a truncation of order
$\delta^{-1/\al}$, which, when $\al\rightarrow 2$, is of the same
order as the one used in obtaining (\ref{f1.25}). However, for the
GUE result, using Gaussian concentration, a truncation of order
$\sqrt{\ln(12b/\delta)}$ gives a slightly better bound, namely,
$$\frac{C_1\sqrt{\ln{\frac{12b}{\delta}}}}{\delta}\exp\bigg\{-\frac{C_2N^2\delta^4}{8ca^2\ln\frac{12b}{\delta}}\bigg\},$$
\noindent where $C_1$ and $C_2$ are absolute constants (different
from those of (\ref{f1.25})).

\end{itemize}
\end{remk}

\vspace{.14in}

Wishart matrices are of interest in many contexts, in particular as
the sample covariance matrix in statistics. Recall that ${\bf
M}={\bf Y}^*{\bf Y}$ is called a complex Wishart matrix if ${\bf Y}$
is a $K \times N $ matrix, $K>N$, with entries ${\bf Y}_{i,j}={\bf
Y}_{i,j}^R+\sqrt{-1}{\bf Y}_{i,j}^I$ (a real Wishart matrix is
defined similarly with ${\bf Y}_{i,j}^I=\delta_0$ and ${\bf M}={\bf
Y}^t{\bf Y}$). Recall also that if the entries of $\bf Y$ are iid
centered random variables with finite variance $\sigma^2$, the
empirical distribution of the eigenvalues of ${\bf Y}^*{\bf Y}/N$
converges as $K\rightarrow\infty$, $N\rightarrow\infty$, and
$K/N\rightarrow\gamma\in(0,+\infty)$ to the
Mar$\check{\text{c}}$enko-Pastur law ([\ref{r38}], [\ref{r29}]) with
density
$$p_\gamma(x)=\frac{1}{2\pi x\gamma\sigma^2}\sqrt{(c_2-x)(x-c_1)},\ \ c_1\le x \le c_2,$$
where $c_1=\sigma^2(1-\gamma^{-1/2})^2$ and
$c_2=\sigma^2(1+\gamma^{-1/2})^2$. When the entries of $\bf Y$ are
iid Gaussian, Johansson [\ref{r27}] and Johnstone [\ref{r30}]
showed, in the complex and real case respectively, that the properly
normalized largest eigenvalue converges in distribution to the
Tracy-Widom law ([\ref{r18}], [\ref{r33}]). Soshnikov [\ref{r34}]
extended the result of Johnstone to Wishart matrix with Non-Gaussian
entries under the condition that $K-N=O(N^{1/3})$ and that the
moments of the entries do not grow too fast. Soshnikov and Fyodorov
[\ref{r17}] recently studied the distribution of the largest
eigenvalue of the Wishart matrix $\bf Y^*Y$, when the entries of
$\bf Y$ are iid Cauchy random variables. We are interested here in
concentration for the linear statistics of the spectral measure and
for the largest eigenvalue of the Wishart matrix $\bf Y^*Y$, where
the entries of $\bf Y$ form an infinitely divisible and, in
particular, a stable one. We restrict our work to the complex
framework, the real framework being essentially the same.

It is not difficult to see that if $\bf Y$ has iid Gaussian entries,
$\bf Y^*Y$ has infinitely divisible entries, each with a L\'evy
measure without a known explicit form. However the dependence
structure among the entries of $\bf Y^*Y$ prevents the vector of
entries to be, itself, infinitely divisible (this is a well known
fact originating with L\'evy, see [\ref{r14}]). The methodology we
previously used cannot be directly applied to deal with functions of
eigenvalues of $\bf Y^*Y$. However, concentration results can be
obtained when we consider the following facts, due to Guionnet and
Zeitouni [\ref{r4}] and already used for that purpose in their
paper.

Let
\begin{equation}\label{f1.20}
A_{i,j} =
\begin{cases}
 0        & \text{for $1\le i\le K,1\le j\le K$} \\
 0        & \text{for $N+1\le i\le K+N, K+1\le j\le K+N$}\\
 1        & \text{for $1\le i\le K, K+1\le j\le K+N$} \\
 1        & \text{for $N+1\le i\le K+N,1\le j\le K$},
 \end{cases}
\end{equation}
and
\begin{equation}
 \omega_{i,j} =
 \begin{cases}
  0        & \text{for $1\le i\le K,1\le j\le K$} \\
  0        & \text{for $N+1\le i\le K+N, K+1\le j\le K+N$}\\
  \bar{Y}_{i,j}       & \text{for $1\le i\le K, K+1\le j\le K+N$} \\
  Y_{i,j}        & \text{for $N+1\le i\le K+N,1\le j\le K$},
  \end{cases}
\end{equation}

\noindent then ${\bf {X_A}}= \left(
  \begin{array}{cc}
    {\bf 0} & {\bf Y}^* \\
    {\bf Y} & {\bf 0} \\
  \end{array}
\right)\in \mathcal M_{(K+N)\times(K+N)}(\mathbb C),$ and
$${{\bf X}^2_{\bf A}}= \left(
   \begin{array}{cc}
     {\bf Y}^*{\bf Y} & {\bf 0} \\
     {\bf 0} & {\bf YY}^* \\
   \end{array}
 \right).$$

\noindent Moreover, since the spectrum of $\bf Y^*Y$ differs from
that of $\bf YY^*$ only by the multiplicity of the zero eigenvalue,
for any function $f$, one has
$$tr(f({{\bf X}^2_{\bf A}}))=2tr(f({\bf Y^*Y}))+(K-N)f(0),$$
and
$$\lambda_{max}({\bf M}^{1/2})=\underset{1\le i\le N}{\max}|\lambda_{i}({\bf {X_A}})|,$$
where ${\bf M}^{1/2}$ is the unique positive semi-definite square
root of ${\bf M}=\bf Y^*Y$.

Next let $\mathbb P^{K,N}$ be the joint law of $({\bf
Y}_{i,j}^R,{\bf Y}_{i,j}^I)_{1\le i\le  K,1\le j\le N}$ on
$\R^{2KN}$, and let $\mathbb E^{K,N}$ be the corresponding
expectation. We present below, in the infinitely divisible case, a
concentration result for the largest eigenvalue $\lambda_{max}({\bf
M})$, of the Wishart matrices ${\bf M}=\bf Y^*Y$. The concentration
for the linear statistic $tr_N(f({\bf M}))$ could also be obtained
using the above observations. \vspace{1in}
\begin{coro}\label{coro1.10}
Let ${\bf M}=\bf Y^*Y$, with ${\bf Y}_{i,j}={\bf
Y}_{i,j}^R+\sqrt{-1}{\bf Y}_{i,j}^I$.

\begin{itemize}
\item[(i)]
Let $X=({\bf Y}_{i,j}^R,{\bf Y}_{i,j}^I)_{1\le i\le N, 1\le j\le K}$
be a random vector with joint law $\mathbb P^{K,N}\sim
ID(\beta,0,\nu)$ such that $\mathbb E^{K,N}[e^{t\|X\|}]<+\infty$,
for some $t>0$. Let
$T=\sup\{t>0:\mathbb{E}^{K,N}[e^{t\|X\|}]<+\infty\}$ and let
$h^{-1}$ be the inverse of $$h(s)=\int_{\R^{2KN}}\|u\| (e^{s
\|u\|}-1)\nu(du),\ \ \ \ 0<s<T.$$ \noindent Then,
\begin{equation}
\mathbb{P}^{K,N}\left(\lambda_{max}({\bf
M}^{1/2})-\mathbb{E}^{K,N}[\lambda_{max}({\bf
M}^{1/2})]\ge\delta\right)\le
e^{-\int_0^{\delta/\sqrt2}h^{-1}(s)ds},
\end{equation}
\noindent for all $0<\delta<h\left(T^{-}\right)$.

\item[(ii)] Let $X=({\bf Y}_{i,j}^R,{\bf Y}_{i,j}^I)_{1\le i\le K,1\le j\le N}$ be
an $\alpha$-stable random vector with L\'evy measure $\nu$ given by
$\nu(B)=\int_{S^{2KN-1}}\sigma(d\xi)\int_0^{+\infty}\mathbf
1_B(r\xi)dr/r^{1+\alpha}$. Then,
$$\mathbb{P}^{K,N}\Bigl(\lambda_{max}({\bf M}^{1/2})-m(\lambda_{max}({\bf M}^{1/2}))\ge\delta\Bigr)\le
C(\al)(\sqrt 2)^\al\frac{\sigma(S^{2KN-1})}{\delta^\al},$$ whenever
$\delta >\sqrt 2a\left[2\sigma(S^{2KN-1})C(\al)\right]^{1/\al}$ and
where $C(\al)=4^\al(2-\al+e\al)/\al(2-\al)$.

\end{itemize}

\end{coro}

\begin{remk}\label{remk1.12}
\begin{itemize}
\item[(i)]
As already mentioned, Soshnikov and Fyodorov ([\ref{r17}]) studied
the asymptotic for the largest singular value of the $K\times N$
random matrix $\bf Y$, which is the largest eigenvalue of the
Wishart matrix $\bf Y^*Y$, when the entries of $\bf Y$ are iid
Cauchy random variables. They argue that although the typical
eigenvalues of $\bf Y^*Y$ is of the order $KN$, the correct order of
the largest eigenvalue of such a matrix is $K^2N^2$. Our result
implies that the largest eigenvalue $\lambda_{max}({\bf M})$ of the
Wishart matrix ${\bf M}=\bf Y^*Y$, when the entries of $\bf Y$ form
an $\al$-stable random vector, is of order at most
$\sigma(S^{2KN-1})^{2/\al}$. We also have a lower bound result which
is described next. In particular, if the entries of the matrix $\bf
Y$ are iid $\al$-stable random variables, the largest eigenvalue of
$\bf Y^*Y$ is of order $K^{2/\al}N^{2/\al}$.

\item[(ii)] Let $X\sim ID(\beta,0,\nu)$ in $\R^d$,
then (see Lemma 5.4 in [\ref{r35}]) for any $x>0$, and any norm
$\|\cdot\|_{\cal N}$ on $\R^d$,
$$\mathbb P\big(\|X\|_{\cal N}\ge x\big)\ge \frac{1}{4}\Big(1-\exp\Big\{-\nu\big(\big\{u\in \R^d: \|u\|_{\cal N}\ge 2x\big\}\big)\Big\}\Big).$$
\noindent But, $\lambda_{max}({\bf M}^{1/2})$ is a norm of the
vector $X=({\bf Y}_{i,j}^R,{\bf Y}_{i,j}^I)$, which we denote by
$\|X\|_\lambda$, if $X$ is a stable vector in $\R^{2KN}$.

\begin{align}
&\mathbb{P}^{K,N}\Bigl(\lambda_{max}({\bf
M}^{1/2})-m(\lambda_{max}({\bf M}^{1/2}))\ge\delta\Bigr)\nonumber\\
&=\mathbb{P}^{K,N}\Bigl(\lambda_{max}({\bf
M}^{1/2})\ge\delta+m(\lambda_{max}({\bf M}^{1/2}))\Bigr)\nonumber\\
&\ge \frac{1}{4}\Big(1-\exp\Big\{-\nu\big(\big\{\lambda_{max}({\bf
M}^{1/2})\ge 2\big(\delta+m(\lambda_{max}({\bf
M}^{1/2}))\big)\big\}\big)\Big\}\Big)\nonumber\\
&\ge\frac{1}{4}\Big(1-\exp\Big\{-\nu\big(\big\{\|X\|_\lambda\ge
2\big(\delta+m(\lambda_{max}({\bf
M}^{1/2}))\big)\big\}\big)\Big\}\Big)\nonumber\\
&=\frac{1}{4}\bigg(1-\exp\bigg\{-\frac{\tilde{\sigma}\big(S_{\|\cdot\|_\lambda}^{2KN-1}\big)}{\al\big(\delta+m(\lambda_{max}({\bf
M}^{1/2}))\big)^\al}\bigg\}\bigg),
\end{align}
\noindent where $S_{\|\cdot\|_\lambda}^{2KN-1}$ is the unit sphere
relative to the norm $\|\cdot\|_\lambda$ and where $\tilde\sigma$ is
the spherical part of the L\'evy measure corresponding to this norm.
Moreover, if the components of $X$ are independent, in which case
the L\'evy measure is supported on the axes of $\R^{2KN}$,
$\tilde{\sigma}\big(S_{\|\cdot\|_\lambda}^{2KN-1}\big)$ is of order
$KN$, and so the largest eigenvalue of ${\bf M}^{1/2}$ is of order
$K^{1/\al}N^{1/\al}$.

\item[(iii)] For any function $f$ such that $g(x)=f(x^2)$ is Lipschitz with
Lipschitz constant $\|g\|_{Lip}:=|||f|||_{\mathcal L}$, $tr(g({\bf
{X_A}}))=tr(f({{\bf X}^2_{\bf A}}))$ is a
 Lipschitz function of the entries of ${\bf Y}$ with Lipschitz constant at most $\sqrt 2|||f|||_{\mathcal
 L}\sqrt{K+N}$. Hence, under the assumptions of part (i) of Corollary
 \ref{coro1.10},
\begin{align}
\mathbb{P}^{K,N}\Big(tr_N(f({\bf M}))&-\mathbb{E}^{K,N}[tr_N(f({\bf
M}))]\ge\delta\frac{K+N}{N}\Big) \nonumber \\
&\le \exp\bigg\{-\int_0^{\sqrt{2(K+N)}\delta\big/|||f|||_{\mathcal
L}}h^{-1}(s)ds\bigg\},
\end{align}
for all $0<\delta<|||f|||_{\mathcal
L}h\left(T^{-}\right)/\sqrt{2(K+N)}$.
\item[(iv)] Under the assumptions of part (ii) of Corollary \ref{coro1.10},
for any function $f$ such that $g(x)=f(x^2)$ is Lipschitz with
$\|g\|_{Lip}=|||f|||_{\mathcal L}$, any median $m(tr_N(f({\bf M})))$
of $tr_N(f({\bf M}))$,
\begin{align}
\mathbb{P}^{K,N}\biggl(tr_N(f({\bf M}&))-m(tr_N(f({\bf M})))\ge\delta\frac{K+N}{N}\biggr) \nonumber \\
&\le C(\al)\frac{|||f|||^\al_{\mathcal
L}}{\sqrt{2^\al(K+N)^\al}}\frac{\sigma(S^{2KN-1})}{\delta^\al},
\end{align}
whenever $\delta>|||f|||_{\mathcal
L}\left[2\sigma(S^{2KN-1})C(\al)\right]^{1/\al}/\sqrt{2(K+N)}$, and
where $C(\al)=4^\al(2-\al+e\al)/\al(2-\al)$.

\end{itemize}

\end{remk}

\begin{remk}
The methodology used to obtain the results of the present paper, in
the absence of the finite exponential moments, can be applied to any
matrices whose entries on and above the main diagonal form such an
infinitely divisible vector $X$. However, to obtain explicit
estimates, we do need specific bounds on $V^2(r)$ and
$\bar{\nu}(r)$, which are not always available when further
knowledge on the L\'evy measure of $X$ is lacking.
\end{remk}

\section{Proofs:}

 We start with a proposition, which is a direct consequence of the
 concentration inequalities obtained in [\ref{r5}] for general
 Lipschitz function of infinitely divisible random vectors with
 finite exponential moment.

 \begin{prop}\label{prop2.1}
 Let $X=(\omega_{i,i}^R,\omega_{i,j}^R,\omega_{i,j}^I)_{1\le i<j\le
 N}$ be a random vector with joint law $\mathbb{P}^N\sim ID(\beta,0,\nu)$
 such that $\mathbb{E}^N \big[e^{t\|X\|}\big]<+\infty$, for some $t>0$ and let
 $T=\sup\{t>0:\mathbb{E}^N \big[e^{t\|X\|}\big]<+\infty\}$. Let $h^{-1}$
 be the inverse of $$h(s)=\int_{\R^{N^2}}\|u\| \big(e^{s\|u\|}-1\big)\nu(du),\ \ \ \ 0<s<T.$$

 \begin{itemize}
 \item[(i)] For any Lipschitz function $f$,
 $$\mathbb{P}^N\!\left(tr_N(f({\bf {X_A}}))-\mathbb{E}^N[tr_N(f({\bf {X_A}}))]\!\ge\!\delta\right)\!\le\! \exp\Bigg\{-\int_0^{\frac{N\delta}{\sqrt2a\flip}}\!h^{-1}(s)ds\Bigg\},$$
 for all $0<\delta< \sqrt 2a\flip h\left(T^{-}\right)/N$.
 \item[(ii)] Let $\lambda_{max}({\bf {X_A}})$ be the largest eigenvalue of the matrix
 ${\bf {X_A}}$. Then,
 $$\mathbb{P}^N\left(\lambda_{max}({\bf {X_A}})-\mathbb{E}^N[\lambda_{max}({\bf {X_A}})]\ge\delta\right)\le \exp\Bigg\{-\int_0^{\frac{\sqrt N\delta}{\sqrt 2a}}h^{-1}(s)ds\Bigg\},$$
 for all $0<\delta<\sqrt 2a h\left(T^{-}\right)/\sqrt{N}$.
 \end{itemize}

 \end{prop}

\vspace{.15in} \noindent{\bf Proof of Theorem \ref{theo1.1}:}

For part (i), following the proof of Theorem 1.3 of [\ref{r4}],
without loss of generality, by shift invariance, assume that
$\min\{x:x\in\mathcal K\}=0$. Next, for any $v>0$, let
\begin{equation}\label{f2.7}
g_{v}(x) =
\begin{cases}
 0        & \text{if $x\le 0$} \\
 x & \text{if $0< x<v$}\\
 v & \text{if $x\ge v$}.
 \end{cases}
\end{equation}
\noindent Clearly $g_v\in Lip(1)$ with $\|g_v\|_\infty=v$. Next for
any function $f\in {Lip_{\mathcal K}(1)}$, any $\Delta>0$, define
recursively $f_\Delta(x)=0$ for $x\le 0$, and for $(j-1)\Delta\le
x\le j\Delta$, $j=1,\dots,\lceil\frac{x}{\Delta}\rceil$, let
$$f_\Delta(x)=\sum_{j=1}^{\lceil\frac{x}{\Delta}\rceil}g_\Delta^{(j)},$$
\noindent where $g_\Delta^{(j)}:=(2\mathbf
1_{\{f(j\Delta)>f_\Delta((j-1)\Delta)\}}-1)g_{\Delta}(x-(j-1)\Delta)$.
Then $|f-f_\Delta|\le \Delta$ and the 1-Lipschitz function
$f_\Delta$ is the sum of at most $|\mathcal K|/\Delta$ functions
$g_\Delta^{(j)}\in Lip(1)$, regardless of the function $f$. Now, for
$\delta>2\Delta$,
\begin{align}
&\mathbb{P}^N \left(\underset{f\in {Lip_{\mathcal K}(1)}}{\sup}\big|tr_N(f({\bf {X_A}}))-\mathbb E^N[tr_N(f({\bf {X_A}}))]\big|\ge\delta\right) \nonumber \\
&\le \mathbb{P}^N \bigg(\underset{f\in {Lip_{\mathcal K}(1)}}{\sup}\bigg\{\big|tr_N(f_\Delta({\bf {X_A}}))-\mathbb E^N(tr_N(f_\Delta({\bf {X_A}})))\big|+\big|tr_N(f({\bf {X_A}}))\nonumber\\
&\ \ \ \ \ \ \ \ \ \ \ \ \ -tr_N(f_\Delta({\bf {X_A}}))\big|+\big|\mathbb E^N[tr_N(f({\bf {X_A}}))]-\mathbb E^N[tr_N(f_\Delta({\bf {X_A}}))]\big|\bigg\}\ge\delta\bigg) \nonumber\\
&\le \mathbb{P}^N
\bigg(\underset{f_\Delta}{\sup}\big|tr_N(f_\Delta({\bf
{X_A}}))-\mathbb E^N(tr_N(f_\Delta({\bf
{X_A}})))\big|>\delta-2\Delta\bigg)\nonumber
\end{align}
\begin{align}\label{f2.27}
&\le \frac{|\mathcal K|}{\Delta}\underset{g^{(j)}_\Delta \in Lip(1)}{\sup}\mathbb{P}^N\left(\big|tr_N(g^{(j)}_\Delta ({\bf {X_A}}))-\mathbb E^N[tr_N(g^{(j)}_\Delta({\bf {X_A}}))]\big|\ge\frac{\Delta(\delta-2\Delta)}{|\mathcal K|}\right) \nonumber\\
&\le \frac{8|\mathcal K|}{\delta}\exp
\bigg\{-\int_0^{\frac{N\delta^2}{8\sqrt2a|\mathcal
K|}}h^{-1}(s)ds\bigg\},
\end{align}
whenever $0<\delta< \sqrt{8\sqrt2a|\mathcal
K|h\left(T^{-}\right)/N}$,
\noindent and where the last inequality
follows from part (i) of the previous proposition by taking also
$\Delta=\delta/4$.

In order to prove part (ii), for any $f\in Lip_b(1)$, i.e, such that
$\flip\le1$, $\|f\|_\infty\le b$, and any $\tau>0$, let $f_{\tau}$
be given via:

\begin{equation}\label{f2.1}
f_\tau(x) =
\begin{cases}
 f(x)        & \text{if $|x|<\tau$} \\
 f(\tau)-\text{sign}(f(\tau))(x-\tau) & \text{if $\tau\le x<\tau+|f(\tau)|$}\\
 f(-\tau)+\text{sign}(f(-\tau))(x+\tau) & \text{if $-\tau-|f(-\tau)|<x\le-\tau$} \\
 0 & \text{otherwise}.
 \end{cases}
\end{equation}
Clearly $f_\tau\in Lip(1)$ and
$supp(f_\tau)\subset\left[-\tau-|f(-\tau)|, \tau+|f(\tau)|\right]$.
\noindent Moreover,
\begin{align}
&\underset{f\in Lip_b(1)}{\sup}\Big|tr_N(f({\bf {X_A}}))\!-\!\mathbb
E^N(tr_N(f({\bf
{X_A}})))\Big|\nonumber\\
&\le \underset{f\in Lip_b(1)}{\sup}\Big|tr_N(f_\tau({\bf
{X_A}}))\!-\!\mathbb E^N(tr_N(f_\tau({\bf
{X_A}})))\Big|\nonumber\\
&\ \ \ \ \ \ \ \ \ \ \ \ \ \ \ +\underset{f\in
Lip_b(1)}{\sup}\Big|tr_N(f({\bf {X_A}})-f_\tau({\bf
{X_A}}))-\mathbb E^N[tr_N(f({\bf {X_A}})-f_\tau({\bf {X_A}}))]\Big| \nonumber \\
&\le \underset{f\in Lip_b(1)}{\sup}\Big|tr_N(f_\tau({\bf
{X_A}}))-\mathbb E^N(tr_N(f_\tau({\bf
{X_A}})))\Big|\nonumber\\
&\ \ \ \ \ \ \ \ \ \ \ \ \ \ \ +2tr_N(g_b(|{\bf
{X_A}}|-\tau))+\!2\mathbb E^N[tr_N(g_b(|{\bf {X_A}}|-\tau))],
\end{align}

\noindent with $g_b$ given as in (\ref{f2.7}). Now,
\begin{align}
&\mathbb P^N\bigg(\underset{f\in Lip_b(1)}{\sup}|tr_N(f({\bf {X_A}}))-\mathbb E^N(tr_N(f({\bf {X_A}})))|\ge\delta\bigg) \nonumber\\
&\le \mathbb P^N\Big(\underset{f\in Lip_b(1)}{\sup}|tr_N(f_\tau({\bf
{X_A}}))-\mathbb E^N(tr_N(f_\tau({\bf {X_A}})))|\ge
\frac{\delta}{3}\Big)\nonumber\\
&\ \ \ \ \ \ \ \ \ \ \ \ \ \ \ \ \ +\mathbb P^N\Bigl(2tr_N(g_b(|{\bf
{X_A}}|-\tau))+2\mathbb E^N[tr_N(g_b(|{\bf
{X_A}}|-\tau))]\ge\frac{2\delta}{3}\Bigr)\nonumber
\end{align}

\begin{align}\label{f2.9}
&\le \mathbb P^N\Big(\underset{f\in Lip_b(1)}{\sup}|tr_N(f_\tau({\bf
{X_A}}))\!-\!\mathbb E^N(tr_N(f_\tau({\bf {X_A}})))|\ge\!
\frac{\delta}{3}\Big)\nonumber\\
&\!+\!\mathbb P^N\!\Bigl(\!tr_N(g_b(|{\bf
{X_A}}\!|\!-\!\tau))\!-\!\mathbb E^N\![tr_N(g_b(|{\bf
{X_A}}\!|\!-\!\tau))]\!\ge\!\frac{\delta}{3}\!-\!2\mathbb
E^N\![tr_N(g_b(|{\bf {X_A}}\!|\!-\!\tau))]\!\Bigr)\nonumber\\
&\le \mathbb P^N\Big(\underset{f\in Lip_b(1)}{\sup}|tr_N(f_\tau({\bf
{X_A}}))\!-\!\mathbb E^N(tr_N(f_\tau({\bf {X_A}})))|\ge\!
\frac{\delta}{3}\Big)\nonumber\\
&\!+\!\mathbb P^N\!\Bigl(\!tr_N(g_b(|{\bf
{X_A}}|\!-\!\tau))\!-\!\mathbb E^N[tr_N(g_b(|{\bf
{X_A}}|\!-\!\tau))]\!\ge\!\frac{\delta}{3}-2b\mathbb
E^N\![tr_N(\mathbf 1_{\{|{\bf {X_A}}|\ge\tau\}}]\Bigr).
\end{align}

\noindent Let us first bound the second probability in (\ref{f2.9}).
Recall that the spectral radius $\rho({\bf {X_A}})=\underset{1\le
i\le N}{\max}|\lambda_i|$ is a Lipschitz function of $X$ with
Lipschitz constant at most $a\sqrt {2/N}$. Hence, for any $0<t\le
T$, and $\gamma>0$ such that $\bar{\nu}(p_\gamma)\le 1/4$,
\begin{align}\label{f2.12}
\mathbb E^N[tr_N(&\mathbf 1_{\{|{\bf
{X_A}}|\ge\tau\}})]=\frac{1}{N}\sum_{i=1}^N\mathbb
P^N\Bigl(|\lambda_i({\bf X_A})|\ge \tau\Bigr)\nonumber\\
&\le \mathbb P^N\big(\rho({\bf {X_A}})\ge \tau\big)\nonumber\\
&\le \mathbb P^N\bigg(\frac{\sqrt N}{\sqrt2a}\rho({\bf {X_A}})-\frac{\sqrt N}{\sqrt2a}\mathbb E^N\big[\rho({\bf {X_A}})\big]\ge \frac{\sqrt N}{\sqrt2a}\tau-G_2(\gamma)\bigg)\nonumber\\
&\le \exp\bigg\{H(t)-\bigg(\frac{\sqrt
N}{\sqrt2a}\tau-G_2(\gamma)\bigg)t\bigg\}
\end{align}
where we have used Proposition \ref{prop1.1} in the next to last
inequality and where the last inequality follows from Theorem 1 in
[\ref{r5}] (p. 1233) with
$$
H(t)=\int_0^th(s)ds=\int_{\R^{N^2}}\big(e^{t\|u\|}-t\|u\|-1\big)\nu(du).
$$
\noindent We want to choose $\tau$, such that $\mathbb
E^N[tr_N(\mathbf 1_{\{|{\bf {X_A}}|\ge\tau\}})]\le \delta/12b$. This
can be achieved if
\begin{equation}\label{f2.26}
\frac{\sqrt N}{\sqrt2a}\tau-G_2(\gamma)\ge
\frac{\ln\frac{12b}{\delta}+H(t)}{t}.
\end{equation}
\noindent Since
$$
\frac{d}{dt}\bigg(\frac{\ln\frac{12b}{\delta}+H(t)}{t}\bigg)=\frac{th(t)-\ln\frac{12b}{\delta}-H(t)}{t^2},
$$
\noindent and
$$
\frac{d^2}{dt^2}\bigg(\frac{\ln\frac{12b}{\delta}+H(t)}{t}\bigg)=\frac{t^3H^{''}(t)-2t(th(t)-\ln\frac{12b}{\delta}-H(t))}{t^4},
 $$
\noindent it is clear that the right hand side of (\ref{f2.26}) is
minimized when $t=t_0$, where $t_0$ is the solution of
$$
th(t)-H(t)-\ln\frac{12b}{\delta}=0,
$$
\noindent and the minimum is then $h(t_0)$.

\noindent Thus, if
\begin{equation}\label{f2.8}
\tau=C_0(\delta,b):=\frac{\sqrt2a}{\sqrt
N}\bigg(G_2(\gamma)+h(t_0)\bigg),
\end{equation}

\noindent then
$$\mathbb E^N[tr_N(\mathbf 1_{\{|{\bf
{X_A}}|\ge\tau\}})]\le \frac{\delta}{12b},$$ \noindent and so,
\begin{align}\label{f2.23}
\mathbb P^N\!\biggl(\!tr_N(g_b(&|{\bf
{X_A}}\!|\!-\!\tau))\!-\!\mathbb E^N\![tr_N(g_b(|{\bf
{X_A}}\!|\!-\!\tau))]\!\ge\!\frac{\delta}{3}\!-\!2b\mathbb
E^N\![tr_N(\mathbf 1_{\{|{\bf {X_A}}|\ge\tau\}}]\biggr)\nonumber\\
&\le \mathbb P^N\!\biggl(\!tr_N(g_b(|{\bf
{X_A}}|\!-\!\tau))\!-\!\mathbb E^N\![tr_N(g_b(|{\bf
{X_A}}|\!-\!\tau))]\!\ge\!\frac{\delta}{6}\biggr) \nonumber\\
&\le\exp\Bigg\{-\int_0^{\frac{N\delta}{6\sqrt2a}}\!h^{-1}(s)ds\Bigg\},
\end{align}
\noindent for all $0<\delta< 6\sqrt 2ah\left(T^{-}\right)/N$, where
Proposition \ref{prop2.1} is used in the last inequality.

For $\tau$ chosen as in (\ref{f2.8}), let $\mathcal
K=[-\tau-b,\tau+b]$, then for any $f\in Lip_b(1)$, $f_\tau \in
{Lip_{\mathcal K}(1)}$. By part (i), the first term in (\ref{f2.9})
is such that
\begin{align}\label{f2.24}
\mathbb P^N\Big(\underset{f\in Lip_b(1)}{\sup}&|tr_N(f_\tau({\bf
{X_A}}))\!-\!\mathbb E^N(tr_N(f_\tau({\bf {X_A}})))|\ge\frac{\delta}{3}\Big)\nonumber\\
&\le \mathbb P^N\bigg(\underset{f_\tau \in {Lip_{\mathcal
K}(1)}}{\sup}|tr_N(f_\tau({\bf {X_A}}))-\mathbb E^N[tr_N(f_\tau({\bf
{X_A}}))]|\ge\frac{\delta}{3}\bigg)\nonumber\\
&\le\frac{48(C_0(\delta,b)+b)}{\delta}\exp\Bigg\{-\int_0^{\frac{N\delta^2}{144\sqrt2a(C_0(\delta,b)+b)}}h^{-1}(s)ds\Bigg\},
\end{align}
\noindent for all $0\!<\!\delta^2\!\le
144\sqrt2a\big(C_0(\delta,b)+b\big)h(T^-)/N$.

Hence, returning to (\ref{f2.9}), using (\ref{f2.23}) and
(\ref{f2.24}) and for $$\delta<\min\Big\{6\sqrt
2ah\left(T^{-}\right)/N,\sqrt{144\sqrt2a\big(C_0(\delta,b)+b\big)h(T^-)/N}\Big\},$$
\noindent we have
\begin{align}
&\mathbb P^N\bigg(\underset{f\in Lip_b(1)}{\sup}|tr_N(f({\bf {X_A}}))-\mathbb E^N(tr_N(f({\bf {X_A}})))|\ge\delta\bigg) \nonumber\\
&\le\!2\frac{24(C_0(\delta,b)\!+\!b)}{\delta}\exp\!\Bigg\{\!\!-\!\!\int_0^{\frac{N\delta}{6\sqrt2a}\frac{\delta}{24(C_0(\delta,b)\!+\!b)}}\!h^{-1}\!(s)ds\Bigg\}\!+\! \exp\!\Bigg\{\!\!-\!\!\int_0^{\frac{N\delta}{6\sqrt2a}}\!h^{-1}\!(s)ds\!\Bigg\}   \nonumber\\
&\le\bigg(2+\frac{1}{12}\bigg)\frac{24(C_0(\delta,b)+b)}{\delta}\exp\Bigg\{-\int_0^{\frac{N\delta^2}{144\sqrt2a(C_0(\delta,b)+b)}}h^{-1}(s)ds\Bigg\},
\end{align}
\noindent since only the case $\delta\le 2b$ presents some interest
(otherwise the probability in the statement of the theorem is zero).
Part (ii) is then proved.

\hfill $\Box$

\vspace{.15in}

\noindent{\bf Proof of Proposition \ref{theo1.3}:}

As a function of $x\in\R^{N^2}$, $d_W(\hat{\mu}_A^N, \mu)(x)$ is
Lipschitz with Lipschitz constant at most $\sqrt 2a/N$. Indeed, for
$x,y\in\R^{N^2}$,

\begin{align}
d_W(\hat{\mu}_A^N, \mu)(x)&=\underset{f \in Lip_b(1)}{\sup}\left|tr_N\big(f({\bf {X_A}})(x)\big)-\int_\R fd\mu\right|  \nonumber\\
&\le \!\!\underset{f\in Lip_b(1)}{\sup}\bigg|tr_N(f({\bf {X_A}})(x))-tr_N(f({\bf {X_A}})(y))\bigg| \nonumber\\
&\ \ \ \ \ \ \ \ \ \ \ \ \ \ \ \ \ \ \ \ \ \ \ \ \ \ \ \ \ \ \ \ +\!\!\underset{f\in Lip_b(1)}{\sup}\left|tr_N(f({\bf {X_A}})(y))-\int_\R fd\mu\right|\nonumber \\
&\le \frac{\sqrt 2a}{N}\|x-y\|+d_W(\hat{\mu}_A^N, \mu)(y).
\end{align}

\noindent Theorem \ref{theo1.3} then follows from Theorem 1 in
[\ref{r5}]. \hfill $\Box$

\vspace{.15in}

\noindent{\bf Proof of Corollary \ref{coro1.4}:}

For L\'evy measures with bounded support, $\mathbb
E^N\big[e^{t\|X\|}\big]<+\infty$, for all $t\ge0$, and moreover
$$h(t)\le
V^2\bigg(\frac{e^{tR}-1}{R}\bigg).$$ \noindent Hence
$$H(t)=\int_0^th(s)ds\le \frac{V^2}{R^2}\big(s^{tR}-1-tR\big),$$
 \noindent and
$$
\exp\bigg\{-\int_0^xh^{-1}(s)ds\bigg\}\le\exp\bigg\{\frac{x}{R}-\bigg(\frac{x}{R}+\frac{V^2}{R^2}\bigg)\ln\bigg(1+\frac{Rx}{V^2}\bigg)\bigg\}.
$$

\noindent Thus, one can take
$$
C(\delta,b)=C\bigg(\frac{\sqrt 2a}{\sqrt
N}\Big(G_2(\gamma)+\frac{V^2}{R}\big(e^{t_0R}-1\big)\Big)+b\bigg),
$$
\noindent where $t_0$ is the solution, in $t$, of
$$\frac{V^2}{R^2}\Big(tRe^{tR}-e^{tR}+1\Big)=\ln\frac{12b}{\delta}.$$
 \noindent Applying Theorem
\ref{theo1.1} (ii) yields the result.

\hfill  $\Box$

\vspace{.15in}

\vspace{.12in} In order to prove Theorem \ref{theo1.7}, we first
need the following lemma, whose proof is essentially as the proof of
Theorem 1 in [\ref{r7}].
\begin{lem}\label{lemma2.4}
Let $X=({\omega_{i,i}^R},\omega_{i,j}^R,\omega_{i,j}^I)_{1\le i<j\le
N}$ be an $\al$-stable vector, $0<\al<2$, with L\'evy measure $\nu$
given by (\ref{f1.6}). For any $x_0, x_1>0$, let
$g_{x_0,x_1}(x)=g_{x_1}(x-x_0)$, where $g_{x_1}(x)$is defined as in
(\ref{f2.7}). Then,
 $$\mathbb{P}^N\biggl(\Big|tr_N(g_{x_0,x_1}({\bf {X_A}}))-\mathbb E^N[tr_N(g_{x_0,x_1}({\bf
{X_A}}))]\Big|\ge\delta\biggr)\le
C(\al)\frac{a^{\al}\sigmaS}{N^\al\delta^\al},$$ whenever
$\delta^{1+\al}>\big(2\sqrt2a\big)^{1+\al}\sigmaS x_1/\al N^{1+\al}$
and where $C(\al)=2^{5\al/2}(2e\al+2-\al)/\al(2-\al)$.

\end{lem}

\vspace{.15in}
\noindent{\bf Proof of Theorem \ref{theo1.7}}

For part (i), first consider $f\in {Lip_{\mathcal K}(1)}$. Using the
same approximation as in Theorem \ref{theo1.1}, any function $f\in
{Lip_{\mathcal K}(1)}$ can be approximated by $f_\Delta$, which is
the sum of at most $|\mathcal K|/\Delta$ functions
$g^{(j)}_\Delta\in Lip(1)$, regardless of the function $f$. Now, and
as before, for $\delta>2\Delta$,

\begin{align}
&\mathbb{P}^N \left(\underset{f\in {Lip_{\mathcal K}(1)}}{\sup}|tr_N(f({\bf {X_A}}))-\mathbb E^N(tr_N(f({\bf {X_A}})))|\ge\delta\right) \nonumber \\
&\le \frac{|\mathcal
K|}{\Delta}\underset{\underset{j=1,\cdots,\lceil\frac{|\mathcal
K|}{\Delta}\rceil}{g^{(j)}_\Delta \in Lip_b(1)}}{\sup}
\mathbb{P}^N\left(\Big|tr_N(g^{(j)}_\Delta ({\bf {X_A}}))-\mathbb
E^N[tr_N(g^{(j)}_\Delta({\bf {X_A}}))]\Big|\ge
\frac{\Delta(\delta-2\Delta)}{|\mathcal K|}\!\right)\nonumber \\
&\le
 \frac{4|\mathcal K|}{\delta}\frac{8^{\al}a^{\al}C_2(\al)\sigmaS|\mathcal K|^\al}{N^\al\delta^{2\al}},
\end{align}
whenever
\begin{equation}\label{f2.3}
\frac{\delta^2}{8|\mathcal K|}>
\frac{2\sqrt2a}{N}\Big(\frac{\sigmaS\delta}{4\al}\Big)^{\frac{1}{1+\al}},
\end{equation}
\noindent and where the last inequality follows from Lemma
\ref{lemma2.4}, taking also $\Delta=\delta/4$.

For any $f\in Lip_b(1)$, and any $\tau>0$, let $f_{\tau}$ be given
as in (\ref{f2.1}). Then, $f_\tau\in {Lip_{\mathcal K}(1)}$, where
$\mathcal K=[-\tau-b,\tau+b]$, and moreover,

\begin{align}\label{f2.10}
&\mathbb P^N\bigg(\underset{f\in Lip_b(1)}{\sup}|tr_N(f({\bf {X_A}}))-\mathbb E^N(tr_N(f({\bf {X_A}})))|\ge\delta\bigg)\nonumber \\
&\le \mathbb P^N\!\Bigl(\!tr_N(g_{\tau,b}(|{\bf
{X_A}}\!|))\!-\!\mathbb E^N\![tr_N(g_{\tau,b}(|{\bf
{X_A}}\!|))]\!\ge\!\frac{\delta}{3}\!-\!2b\mathbb E^N\![tr_N(\mathbf
1_{\{|{\bf {X_A}}|\ge\tau\}}]\!\Bigr)\nonumber\\
&\ \ \ \ \ \ \ \ \ \ \ \ +\mathbb P^N\Big(\underset{f_\tau \in
{Lip_{\mathcal K}(1)}}{\sup}|tr_N(f_\tau({\bf {X_A}}))-\mathbb
E^N(tr_N(f_\tau({\bf {X_A}})))|\ge \frac{\delta}{3}\Big).
\end{align}

\noindent The spectral radius $\rho({\bf {X_A}})$ is a Lipschitz
function of $X$ with Lipschitz constant at most $\sqrt2a/\sqrt N$.
Then by Theorem 1 in [\ref{r7}],
\begin{align}
\mathbb E^N[tr_N(\mathbf 1_{\{|{\bf
{X_A}}|\ge\tau\}})]&=\frac{1}{N}\sum_{i=1}^N{\mathbb
P^N\Bigl(|\lambda_i({\bf {X_A}})|\ge\tau\Bigr)}\nonumber\\
&\le \mathbb P^N\Big(\rho({\bf {X_A}})>\tau\Big)\nonumber\\
&\le \mathbb P^N\Bigg(\rho({\bf {X_A}})-m(\rho({\bf {X_A}}))>\tau-\frac{\sqrt2a}{\sqrt N}J_1(\al)\Bigg)\nonumber\\
&\le
\frac{C_1(\al)2^{\al/2}a^\al\sigmaS}{N^{\al/2}\big(\tau-\frac{\sqrt2a}{\sqrt
N}J_1(\al)\big)^\al},
\end{align}
whenever
\begin{equation}\label{f2.4}
\bigg(\tau-\frac{\sqrt2a}{\sqrt
N}J_1(\al)\bigg)^\al\ge\frac{2C_1(\al)2^{\al/2}a^\al\sigmaS}{N^{\al/2}},
\end{equation}
\noindent and where $C_1(\al)=4^\al(2-\al+e\al)/\al(2-\al)$.
\noindent Now, if $\tau$ is chosen such that
$$\frac{C_1(\al)2^{\al/2}a^\al\sigmaS}{N^{\al/2}\big(\tau-\frac{\sqrt2a}{\sqrt N}J_1(\al)\big)^\al}\le \frac{\delta}{12b},$$

\noindent that is, if
\begin{equation}\label{f2.5}
\bigg(\tau-\frac{\sqrt2a}{\sqrt
N}J_1(\al)\bigg)^\al\ge\frac{12bC_1(\al)2^{\al/2}a^\al\sigmaS}{\delta
N^{\al/2}},
\end{equation}
\noindent it then follows that $$\mathbb E^N[tr_N(\mathbf 1_{\{|{\bf
{X_A}}|\ge\tau\}})]\le \frac{\delta}{12b}.$$

\noindent Since $g_{\tau,b}(|{\bf {X_A}}\!|)$ is the sum of two
functions of the type studied in Lemma \ref{lemma2.4} with $x_1=b$,
we have,

\begin{align}\label{f2.11}
\mathbb P^N\!\Bigl(tr_N(g_{\tau,b}(&|{\bf {X_A}}|))\!-\!\mathbb
E^N[tr_N(g_{\tau,b}(|{\bf
{X_A}}|))]\!\ge\!\frac{\delta}{3}\!-\!2b\mathbb E^N[tr_N(\mathbf
1_{\{|{\bf {X_A}}|\ge\tau\}}]\Bigr)\nonumber\\
&\le 2\mathbb P^N\Bigl(tr_N(g_{\tau,b}({\bf {X_A}}))-\mathbb
E^N[tr_N(g_{\tau,b}({\bf {X_A}}))]\ge\frac{\delta}{12}\Bigr)\nonumber\\
&\le 2C_2(\al)\frac{12^{\al}a^\al\sigmaS}{N^\al\delta^\al},
\end{align}
\noindent whenever
\begin{equation}\label{f2.6}
\delta^{1+\al}>\Big(\frac{2\sqrt2a}{N}\Big)^{1+\al}
\frac{12^{1+\al}\sigmaS b}{\al},
\end{equation}
\noindent and where $C_2(\al)=2^{5\al/2}(2e\al+2-\al)/\al(2-\al)$.
\noindent The respective range (\ref{f2.5}) and (\ref{f2.6}) suggest
that one can choose, for example,
$$\tau=\frac{\sqrt2a}{\sqrt
N}J_1(\al)+\frac{\sqrt2a}{\sqrt N}\delta.$$ \noindent Then, there
exists $\delta(\al,a,N,\nu)$ such that for
$\delta>\delta(\al,a,N,\nu)$,
\begin{equation}
\begin{split}
&\mathbb P^N\Big(\underset{f\in Lip_b(1)}{\sup}|tr_N(f({\bf {X_A}}))-\mathbb E^N[tr_N(f({\bf {X_A}}))]|\ge\delta\Big) \\
&\le \mathbb P^N\Big(\underset{f_\tau {Lip_{\mathcal
K}(1)}}{\sup}|tr_N(f_\tau({\bf {X_A}}))-\mathbb E^N[tr_N(f_\tau({\bf
{X_A}}))]|\ge
\frac{\delta}{3}\Big)\\
&\ \ \ \ \ \ +\mathbb P^N\!\Bigl(tr_N(g_{\tau,b}(|{\bf
{X_A}}|))\!-\!\mathbb E^N[tr_N(g_{\tau,b}(|{\bf
{X_A}}|))]\!\ge\!\frac{\delta}{3}\!-\!2b\mathbb E^N[tr_N(\mathbf
1_{\{|{\bf {X_A}}|\ge\tau\}}]\Bigr)\\
&\le \frac{C_3(\al)a^{\al}\sigmaS\Big(\frac{\sqrt2a}{\sqrt
N}J_1(\al)+b+\frac{\sqrt2a}{\sqrt
N}\delta\Big)^{1+\al}}{N^\al\delta^{1+2\al}}+\frac{C_4(\al)a^{\al}\sigmaS}{N^\al\delta^{\al}},\nonumber
\end{split}
\end{equation}
where $C_3(\al)=2^{4+2\al}12^\al C_2(\al)$,
$C_4(\al)=2(12^{\al})C_2(\al)$ and $\delta(\al,a,N,\nu)$ is such
that (\ref{f2.3}) and (\ref{f2.6}) hold.

Part (ii) is a direct consequence of Theorem 1 of [\ref{r7}], since
$d_W(\hat{\mu}_A^N, \mu)\in Lip(\sqrt2a/N)$ as shown in the proof of
Proposition \ref{theo1.3}. \hfill $\Box$

\vspace{.15in} \noindent{\bf Proof of Theorem \ref{theo1.8}}

For any $f\in Lip(1)$, Theorem 1 in [\ref{r7}] gives a concentration
inequality for $f(X)$, when it deviates from one of its medians. For
$1<\al<2$, a completely similar (even simpler) argument gives the
following result,
\begin{align}\label{f2.53}
\mathbb{P}^N\!\bigl(f(X)-\mathbb E^N[f(X)]\!\ge\! x\bigr)&\le
\frac{C(\al)\sigmaS}{x^\al},
\end{align}
whenever $x^\al\ge\!K(\al)\sigmaS$, where
$C(\al)=2^\al(e\al+2-\al)/(\al(2-\al))$ and
$K(\al)=\max\big\{2^\al/(\al-1),C(\al)\big\}$.

Next, following the proof of Theorem \ref{theo1.1}, approximate any
function $f\in Lip_b(1)$ by $f_\tau\in Lip_{[-\tau-b,\tau+b]}(1)$
defined via (\ref{f2.1}). Hence,

\begin{align} \label{f2.2}
&\mathbb P^N\Big(\underset{f\in Lip_b(1)}{\sup}|tr_N(f({\bf {X_A}}))-\mathbb E^N(tr_N(f({\bf {X_A}})))|\ge\delta\Big) \nonumber\\
&\le \mathbb P^N\Big(\underset{f_\tau \in {Lip_{\mathcal
K}(1)}}{\sup}|tr_N(f_\tau({\bf {X_A}}))-\mathbb E^N(tr_N(f_\tau({\bf
{X_A}})))|\ge
\frac{\delta}{3}\Big)\nonumber\\
&\!+\!\mathbb P^N\!\Bigl(tr_N(g_{\tau,b}(|{\bf {X_A}}|))\!-\!\mathbb
E^N[tr_N(g_{\tau,b}(|{\bf
{X_A}}|))]\!\ge\!\frac{\delta}{3}\!-\!2b\mathbb E^N[tr_N(\mathbf
1_{\{|{\bf {X_A}}|\ge\tau\}}]\Bigr).
\end{align}

\noindent For $\rho({\bf {X_A}})$ the spectral radius of the matrix
${\bf {X_A}}$, and for any $\tau$, such that $\tau-\mathbb
E^N[\rho({\bf {X_A}})]\ge \Big(\frac{\sqrt2a}{\sqrt
N}K(\al)\sigmaS\Big)^{1/\al}$,
\begin{align}\label{f2.54}
\mathbb{E}^N\bigl(tr_N(\mathbf 1_{\{|{\bf {X_A}}|>\tau\}})\bigr)&\le
\mathbb{P}^N\Bigl(\rho({\bf {X_A}})-\mathbb E^N[\rho({\bf
{X_A}})]\ge \tau-\mathbb E^N[\rho({\bf {X_A}})]\Bigr)\nonumber \\
&\le \frac{\Big(\frac{\sqrt2a}{\sqrt N}\Big)^\al
C(\al)\sigmaS}{\big(\tau-\mathbb E^N[\rho({\bf {X_A}})]\big)^\al},
\end{align}
\noindent where we have used, in the last inequality, (\ref{f2.53})
and the fact that $\rho({\bf {X_A}})\in Lip(\frac{\sqrt2a}{\sqrt
N})$. For $Q>0$, let $\tau=\mathbb E^N[\rho({\bf
{X_A}})]+Q\delta^{-1/\al}$. With this choice, we then have:

\begin{align}\label{f2.54}
\mathbb{E}^N\bigl(tr_N(\mathbf 1_{\{|{\bf {X_A}}|>\tau\}})\bigr)
&\le \frac{\Big(\frac{\sqrt2a}{\sqrt N}\Big)^\al
C(\al)\sigmaS}{\big(\tau-\mathbb E^N[\rho({\bf
{X_A}})]\big)^\al}\nonumber\\
&\le \delta\frac{\Big(\frac{\sqrt2a}{\sqrt N}\Big)^\al
C(\al)\sigmaS}{Q^\al}
\nonumber\\
&\le \frac{\delta}{12b},
\end{align}

\noindent provided $Q^\al/\delta>\sqrt2aK(\al)\sigmaS/\sqrt N$, and
$\Big(\frac{\sqrt2a}{\sqrt N}\Big)^\al C(\al)\sigmaS/Q^\al\le
1/(12b)$. Now, taking
$Q=\sqrt2a\big(12bC(\al)\sigmaS\big)^{1/\al}/\sqrt N$, and
recalling, for $1<\al<2$, the lower range concentration result for
stable vectors (Theorem 1 and Remark 3 in [\ref{r2}]): For any
$\epsilon>0$, there exists $\eta_0(\epsilon)$, such that for all
$0<\delta<\sqrt 2a\flip\eta_0(\epsilon)/N$,
\begin{align}
\mathbb{P}^N\bigl(tr_N&(f({\bf {X_A}}))-\mathbb E^N(tr_N(f({\bf
{X_A}})))\ge\delta\bigr)\nonumber\\
 &\le
 (1+\epsilon)\exp\Bigg\{-\frac{\frac{2-\al}{10}\big(\frac{\al-1}{\al}\big)
 ^{\frac{\al}{\al-1}}}{(\sigma(S^{N^2-1}))^{1/(\al-1)}}\Bigg(\frac{N}{\sqrt2a\flip
 }\Bigg)^{\frac{\al}{\al-1}}\delta^{\frac{\al}{\al-1}}\Bigg\}.
\end{align}

\noindent With arguments as in the proof of Theorem \ref{theo1.1},
if
$$\delta<\eta(\epsilon):=\Bigg(\frac{72\sqrt 2a}{N}\bigg(\frac{\sqrt2a}{\sqrt
N}J_2(\al)+b+\bigg(\frac{\sqrt2a}{\sqrt
N}K(\al)\sigmaS\bigg)^{1/\al}\bigg)\eta_0(\epsilon)\Bigg)^{1/2},$$
\noindent there exist constants $D_1\big(\al,a,N,\sigmaS\big)$ and
$D_2\big(\al,a,N,\sigmaS\big)$, such that the first term in
(\ref{f2.2}) is bounded above by

\begin{equation}\label{f2.60}
(1+\epsilon)\frac{D_1\big(\al,a,N,\sigmaS\big)}{\delta^{\frac{\al+1}{\al}}}\exp\Bigl(-D_2\big(\al,a,N,\sigmaS\big)
\delta^{\frac{2\al+1}{\al-1}}\Bigr).
\end{equation}
\noindent Indeed, with the choice of $\tau$ above and $D^*$ as in
(\ref{f1.26}), $2(\tau+b)\le D^*/\delta^{1/\al}$. Moreover, as in
obtaining (\ref{f2.27}), $D_1$ can be chosen to be $24D^*$, while
$D_2$ can be chosen to be
$$\frac{\frac{2-\al}{10}\big(\frac{\al-1}{\al}\big)
 ^{\frac{\al}{\al-1}}}{\Big(\sigmaS\Big)^{\frac{1}{\al-1}}}\bigg(\frac{N}{\sqrt2a}\bigg)^{\frac{\al}{\al-1}}
\frac{1}{\big(72D^*\big)^{\frac{\al}{\al-1}}}.$$

We remind the reader that, as already mentioned, $J_2(\al)$ can be
replaced by $\mathbb E^N[\|X\|]$. According to the result of Marcus
and Rosi\'nski [\ref{r37}] and the estimate in [\ref{r6}], if
$\mathbb E^N[X]=0$, then
$$
\frac{1}{4(2-\al)^{1/\al}}\sigmaS^{1/\al}\le \mathbb E^N[\|X\|]\le
\frac{17}{8\big((2-\al)(\al-1)\big)^{1/\al}}\sigmaS^{1/\al}.
$$
Finally, note that, as in the proof of Theorem \ref{theo1.1} (ii),
the second term in (\ref{f2.2}) is dominated by the first term. The
theorem is then proved, with the constant $D_1\big(a,N,\sigmaS\big)$
magnified by 2.

\hfill $\Box$

\vspace{.15in}

 \noindent{\bf Proof of Corollary \ref{coro1.10}:}

As a function of $({\bf Y}^R_{i,j},{\bf Y}^I_{i,j})_{1\le i\le K,
1\le j\le N}$, with the choice of $A$ made in (\ref{f1.20}),
$\lambda_{max}({\bf {X_A}})\in Lip(\sqrt2)$. Hence part(i) is a
direct application of Theorem 1 in [\ref{r5}], while part(ii) can be
obtained by applying Theorem \ref{theo1.5}. \hfill $\Box$

\vspace{.15in}

{\bf Acknowledgements} Both authors would like to thank the
organizers of the Special Program on High-Dimensional Inference and
Random Matrices at SAMSI. Their hospitality and support, through the
grant DMS-0112069, greatly facilitated the completion of this paper.

 \setcounter{equation}{0}

\end{document}